\documentclass[10pt]{article}  
\usepackage[all]{xy}
\usepackage[mathscr]{euscript}
\usepackage{amsmath, amsfonts, amssymb, mathrsfs}
\usepackage{appendix}
\usepackage{amsthm}
\usepackage{comment}
\usepackage{textcomp}
\usepackage{setspace}
\usepackage{fix-cm}

\usepackage[alphabetic,lite]{amsrefs}

\usepackage{tikz-cd}

\usetikzlibrary{arrows,positioning,decorations.pathmorphing,
  decorations.markings
 }

\tikzset{inner sep=0pt, 
  root/.style={circle,draw,minimum size=7pt,thick}, 
  fatroot/.style={circle,draw,minimum size=10pt,thick}, 
  short root/.style={circle,fill,minimum size=7pt}, 
  doublearrow/.style={postaction={decorate}, 
  decoration={markings,mark=at position .7
  with {\arrow{angle 60}}},double distance=3pt,thick}
}

\newcommand{\ZZ}{\mathbb Z}

\newcommand{\bQ}{\mathbb Q}
\newcommand{\QQ}{\mathbb Q}
\newcommand{\RR}{\mathbb R}

\newcommand{\frakg}{\mathfrak g}
\newcommand{\frakh}{\mathfrak h}
\newcommand{\fraksl}{\mathfrak sl}
\newcommand{\frakgl}{\mathfrak gl}

\newcommand{\frakt}{\mathfrak t}

\newcommand{\la}{\langle}
\newcommand{\ra}{\rangle}

\newcommand{\al}{\alpha}
\newcommand{\be}{\beta}

\newcommand{\lam}{\lambda}
\newcommand{\varep}{\varepsilon}

\DeclareMathOperator{\Ad}{Ad}
\DeclareMathOperator{\Aut}{Aut}

\DeclareMathOperator{\Gal}{Gal}
\DeclareMathOperator{\Hom}{Hom}

\DeclareMathOperator{\GL}{GL}
\DeclareMathOperator{\PGL}{PGL}

\DeclareMathOperator{\Ind}{Ind}

\DeclareMathOperator{\Pic}{Pic}
\DeclareMathOperator{\Div}{Div}

\theoremstyle{plain}
\input xy
\xyoption{all}
\newtheorem{theorem}{Theorem}[section]

\newtheorem{Lem}[theorem]{Lemma}
\newtheorem{Prop}[theorem]{Proposition}

\theoremstyle{remark}
\newtheorem{remark}[theorem]{Remark}
\numberwithin{equation}{section}
\numberwithin{paragraph}{section}

\newcommand{\wt}{\widetilde}
\newcommand{\lieaut}{\wt w} 
\newcommand{\Heis}{\mathscr{H}}

\newcommand{\affS}{\mathscr{S}}
\newcommand{\affC}{\mathscr{C}}

\AtEndDocument{\bigskip{\footnotesize 
\textsc{Beth Romano}~~ \texttt{beth.romano@maths.ox.ac.uk}\\\textsc{Mathematical Institute, University of Oxford,
Andrew Wiles Building,
Woodstock Road,
Oxford,
OX2 6GG}}}

\begin{document}

\author{Beth Romano}
\title{On central extensions and simply laced Lie algebras}
\date{}
\maketitle

\begin{abstract}
Let $\Lambda$ be a simply laced root lattice and $w$ an elliptic automorphism of $\Lambda$ of order $d$. This paper gives a construction that begins with a central extension of the group of coinvariants $\Lambda_w$ and produces a semisimple Lie algebra of Dykin type $\Lambda$ with an automorphism of order $d$ lifting $w$. The input for this construction naturally arises when considering certain families of algebraic curves. The construction generalizes one used by Thorne  to study plane quartics and one used by the author and Thorne to study a family of genus-2 curves.
\end{abstract}

\thispagestyle{empty}  

{
\renewcommand{\thefootnote}{}  
	\footnotetext{{MSC: Primary 17B40, 17B20, 17B10; Secondary 14H40}}
	\footnotetext{{Keywords: Semisimple Lie algebras, Vinberg theory}}
}

\tableofcontents

\section{Introduction}

Graded Lie algebras are a key source of coregular representations and thus have become important for recent developments in arithmetic statistics (see, e.g., \cite{Gross}, \cite{Tho13}, \cite{Tho16}, \cite{RomanoThorne1}, \cite{RomanoThorne2}). In this paper, I give a general construction that begins with a central extension of a certain form and produces a graded Lie algebra. This construction generalizes those of \cite{Tho16} and \cite{RomanoThorne2}. In both \cite{Tho16} and \cite{RomanoThorne2}, the construction was applied to appropriate Mumford-theta groups and allowed us associate orbits in graded Lie algebras to points in quotients of the Jacobians of certain algebraic curves. The present paper, while mainly algebraic in content, develops the generalized contruction as a tool for further number-theoretic applications.

To explain the results in more detail, we let $\Lambda$ be a simply laced root lattice, $w$ an elliptic automorphism of $\Lambda$ of order $d$, and $k$ a field of characteristic zero containing a primitive $d$th root of unity $\zeta$. Let
\begin{equation}\label{eqn-centralext}
1 \to \la \zeta \ra \to \Heis \to \Lambda_w \to 1
\end{equation} 
be a central extension of the coinvariants of $w$, and let $\la \cdot, \cdot \ra$ be the commutator pairing of $\Heis$. Given a bilinear pairing $\varep: \Lambda \times \Lambda \to k^\times$ that satisfies certain properties with respect to the pair $(w, \Heis)$ (see Section \ref{section-input}), I give a construction of 
a Lie algebra $\frakh$ over $k$ as well as an order-$d$ automorphism $\lieaut$ of $\frakh$ lifting $w$. The decomposition of $\frakh$ into eigenspaces for $\lieaut$ gives $\frakh$ the structure of a graded Lie algebra. As in \cite{RomanoThorne2}, the details of the construction of $\frakh$ are inspired by ideas from \cite{Lurie} and \cite{Lepowsky}. 

As an example of the input, for a fixed pair $(\Lambda, w)$ as above, consider the pairing $\la \cdot, \cdot \ra_w : \Lambda \times \Lambda \to \la \zeta \ra$ given by
\begin{equation*}
\la \al, \be \ra_w = \zeta^{\sum_{j = 1}^{d - 1} j(w^j\al, \be)}
\end{equation*}
for all $\al, \be \in \Lambda$, which appears in \cite{Lepowsky}. Then $\la \al, \be \ra_w$ gives a well-defined pairing on $\Lambda_w$, and there exists a central extension of the form (\ref{eqn-centralext}) with commutator pairing $\la \cdot, \cdot \ra_w$. Then the pairing $\varep_w$ defined by
\begin{equation*}
\varep_w(\al, \be) = \prod_{j = 1}^{d - 1} (1 - \zeta^{-j})^{(w^j\al, \be)}
\end{equation*}
satisfies the properties required for the construction (see Section \ref{section-pairingsandexts}). We note that in examples coming from algebraic geometry, the pairing $\la \cdot, \cdot \ra_w$ naturally arises as the Weil pairing on the $d$-torsion points of the Jacobians of certain curves (see \cite{RomanoThorne2}, \cite{Kulkarni}, and Section \ref{section-examples} below).

Further, in the case when the Lie algebra $\frakh$ is constructed using the pairing $\varep_w$, there is a relationship between the representation theory of $\Heis$ and that of the fixed-point subalgebra $\mathfrak{g} = \mathfrak{h}^{\lieaut}$: given a representation $(\rho, V)$ of $\Heis$ over $k$ such that $\rho(\zeta) v = \zeta v$ for all $v \in V$, there is a naturally associated representation of $\frakg$ on $V$ (see Theorem \ref{thm-rep-ext}). The existence of such a representation in the case when $\Lambda$ is type $E_8$ and $d = 3$ is an essential component of the proof of \cite[Theorem 3.7]{RomanoThorne2}, in which we show that $\Heis$ naturally embeds as a subgroup of a simply connected group with Lie algebra $\frakg$.

\subsection{Structure of the paper}
The paper has the following structure. Section \ref{section-construction} gives the details of the construction of $\frakh$ and its grading: in Section \ref{section-input}, we define an input datum as a quadruple $(\Lambda, w, \Heis, \varep)$ of objects as described above. In Section \ref{section-hconstruction}, given an input datum, 
we construct a Lie algebra $\frakh$ with Dynkin type given by $\Lambda$. In Section \ref{section-functor}, 
we show that the construction of $\frakh$ is functorial.
Using this, in Section \ref{section-lieaut} we construct an order-$d$ automorphism $\lieaut$ of $\frakh$ lifting $w$.
In Section \ref{section-galois}, we assume $k/k_0$ is a finite Galois extension
and that $\Gal(k/k_0)$ acts on 
the input datum for $\frakh$ via automorphisms. Assuming this action has certain properties, we show that $\frakh$ and $\frakg$ descend to Lie algebras over $k_0$.

In Section \ref{section-pairingsandexts}, we start with a pair $(\Lambda, w)$ and explicitly construct an input datum $(\Lambda, w, \Heis, \varep_w)$, where $\Heis$ has commutator pairing $\la \cdot, \cdot \ra_w$ and $\varep_w$ is as defined above.
In Section \ref{section-reps-of-g}, we show that with input datum as defined in Section \ref{section-pairingsandexts}, given a representation $(\rho, V)$ of $\Heis$ over $k$ such that $\rho(\zeta)v = \zeta v $ for all $v \in V$, there is an associated representation of $\frakg$ on $V$.

Section \ref{section-lie} contains some Lie-theoretic applications and examples. In particular, we show that the construction of Section \ref{section-construction} gives an alternate proof that given a simple, simply laced Lie algebra $\mathfrak{s}$ with Cartan subalgebra $\frakt$ over an algebraically closed field, every elliptic automorphism of $\frakt$ lifts to an automorphism of $\mathfrak{s}$ of the same order. (This result was already proven in \cite{AdamsHe} using different methods.) 
We also show that using the functoriality described in Section \ref{section-functor}, we may use Section \ref{section-construction} to construct simple Lie algebras of every Dynkin type. As an example of the Galois descent of Section \ref{section-galois}, we show that we may use the construction to obtain a Lie algebra of type $G_2$ over $\bQ$. 

In Section 
\ref{section-examples}, I show that the input for the construction of Section \ref{section-construction} naturally arises when considering certain families of algebraic curves. These families have been studied previously (e.g. in \cite{BS5}, \cite{RainsSam}, \cite{RomanoThorne2}, \cite{Kulkarni}), but we put the relevant ideas into the setting of Section \ref{section-construction}.
We start with 
a semiuniversal deformation $\affS \to B$ of a simple singularity of type $E_8$. For smooth fibers $\affS_\lam$ of this deformation, $\Pic(\affS_\lam)$ is a root lattice of type $E_8$. There are automorphisms of $\affS$ of orders $2, 3$, and $5$ corresponding to elliptic automorphisms of the $E_8$ root lattice. Given a smooth projective curve $C_\lam$ coming from the fixed points of the automorphism of order $d$, the pairing $\la \cdot, \cdot \ra_w$ corresponds to the Weil pairing on $J_{C_\lam}[d]$ where $J_{C_\lam}$ is the Jacobian of $C_\lam$. We may then define central extensions of $J_{C_\lam}[d]$ that give input datum as defined in Section \ref{section-input}.

\subsection{Conventions and notation}

In this paper, we consider a {lattice} $\Lambda$ to be a finitely generated free $\ZZ$-module with a symmetric positive-definite bilinear form $(\cdot, \cdot): \Lambda \times \Lambda \to \ZZ$. We say such a lattice $\Lambda$ is a {simply laced root lattice} if $\Lambda$ is generated as a $\ZZ$-module by $\{\al \in \Lambda \mid (\al, \al) = 2\}$. If $\Lambda$ is a simply laced root lattice, then $\Phi := \{\al \in \Lambda \mid (\al, \al) = 2\}$ forms a simply laced reduced root system in the Euclidean space $\Lambda \otimes_\ZZ \RR$, and we call elements of $\Phi$ {roots}. If $\Lambda$ and $\Lambda'$ are  lattices with corresponding bilinear forms $(\cdot, \cdot)$ and $(\cdot, \cdot)'$ respectively, then an {isomorphism} $\phi: \Lambda \to \Lambda'$ is a $\ZZ$-linear bijection such that $(\phi(\al), \phi(\be))' = (\al, \be)$ for all $\al, \be \in \Lambda$. Two simply laced root lattices are said to be {isomorphic} if they are isomorphic as lattices. 
If $\Lambda$ is simply laced root lattice, then given $\gamma \in \Lambda$, we let $\gamma^\vee \in \Lambda^{\vee} = \Hom(\Lambda, \ZZ)$ be defined by $\gamma^\vee(\al) = (\gamma, \al)$ for all $\al \in \Lambda$. 

\subsection*{Acknowledgements} 
The author received support from EPSRC First Grant EP/N007204/1. This project has received funding
from the European Research Council (ERC) under the European Union's Horizon
2020 research and innovation programme (grant agreement No. 714405).
The author would like to thank Jack Thorne for suggesting the topic, as well as for insightful conversations throughout the writing of the paper. I'd also like to thank Mark Reeder for his helpful comments about an earlier draft of the paper, as well as Jef Laga for helpful conversations about Section \ref{section-examples} and for carefully reading a draft of the paper. I'd like to thank the anonymous referee for their careful reading and valuable suggestions.

\section{Construction of a graded Lie algebra}\label{section-construction}

\subsection{Input data}\label{section-input}

Let $\Lambda$ be a simply laced root lattice with associated bilinear form $(\cdot, \cdot)$ and $\Phi \subset \Lambda$ the subset of roots.
We assume that $(\cdot, \cdot)$ is normalized such that $(\al, \al) = 2$ for all roots $\al$. 
Let $w$ be an elliptic automorphism of $\Lambda$ of order $d$. (Recall 
that $w$ is elliptic if and only if
the group of coinvariants $\Lambda_w := \Lambda/(1-w)\Lambda$ is finite.) 
Let $k$ be a field of characteristic zero, and assume that $k$ contains a primitive $d$th root of unity $\zeta$. 
Let 
\begin{equation*}
1 \to \la \zeta \ra \to \Heis \to \Lambda_w \to 1
\end{equation*}
be a central extension, and let $\la \cdot, \cdot \ra: \Lambda_w \times \Lambda_w \to \la \zeta \ra$ be its commutator pairing. Let $\wt\Lambda = \Lambda \times_{\Lambda_w} \Heis$. Then $\wt\Lambda$ is a central extension
\begin{equation*}
1 \to \la \zeta \ra \to \wt\Lambda \to \Lambda \to 1
\end{equation*}
with commutator pairing given by 
\begin{equation*}
\Lambda \times \Lambda \to \Lambda_w \times \Lambda_w \to \la \zeta \ra
\end{equation*}
where the first map is the natural projection and the second is the pairing $\la \cdot, \cdot \ra$. We continue to write $\la \cdot, \cdot \ra$ for the commutator pairing on $\Lambda$ without risk of confusion. 
Write $\wt\Phi$ for the preimage of $\Phi$ in $\wt\Lambda$ (thus $\wt\Phi$ is a $d$-fold cover of $\Phi$). Given $\wt\al \in \wt\Phi$, we write $\al$ for the image of $\wt\al$ in $\Phi$. Thus by definition, we have
\begin{equation}
\wt\al\wt\be = \la\al, \be\ra\wt\be\wt \al
\end{equation}
for all $\wt\al, \wt\be \in \wt\Phi$.
We note for easy reference that because $\la \cdot, \cdot \ra$ is a commutator pairing, it satisfies the following properties. 

\begin{Lem}\label{lem-comm-props}
The pairing $\la \cdot, \cdot \ra$ has the following properties:
\begin{enumerate} 
\item[1.] $\la \cdot, \cdot \ra$ is bilinear in the sense that for all $\al, \be, \gamma \in \Lambda$,
\begin{itemize}
\item[] 
$\la \al, \be + \gamma\ra = \la \al, \be \ra\la \al, \gamma\ra$,
\item[] 
$\la \al + \be, \gamma \ra = \la \al, \gamma\ra\la \be, \gamma \ra$, and
\item[] 
$\la c\al, \be \ra = \la \al, c\be \ra = \la\al, \be\ra^c$ for all integers $c$.
\end{itemize}
\item[2.] $\la \cdot, \cdot \ra$ is skew symmetric in the sense that $\langle \be, \al \rangle = \la \al, \be \ra^{-1}$ for all $\al, \be \in \Lambda$.
\item[3.]
$\wt\al\wt\be = \wt\be\wt\al$ for all $\wt\al, \wt\be \in \wt\Phi$ such that  $\al + \be = 0$.
\end{enumerate}
\end{Lem}

Let $\varep: \Lambda \times \Lambda \to k^\times$ be a bilinear map (in the sense of Lemma \ref{lem-comm-props}). We say that the  quadruple $(\Lambda, w, \Heis, \varep)$ is an \textit{input datum} if it satisfies the following two properties:
\begin{enumerate}
\item[1.]
$\frac{\varep(\al, \be)}{\varep(\be, \al)} = -\la\be, \al\ra$
whenever $(\al, \be) = -1$. 
\item[2.] $\varep(w\al, w\be) = \varep(\al, \be)$.
\end{enumerate}

In Section \ref{section-pairing}, we describe an alternating pairing $\la \cdot, \cdot \ra_w$ on $\Lambda_w$ previously defined in \cite{Lepowsky} and \cite{ReederEll}. We then describe
one choice for the pairing $\varep$ under the assumption that $\la \cdot, \cdot \ra_w$ is the commutator pairing on $\Heis$.
We postpone this discussion because in many cases, for a fixed commutator pairing $\la \cdot, \cdot \ra$, more than one pairing $\varep$ will satisfy the conditions required. For example, in \cite{RomanoThorne2}, we have $d = 3$ and take $\varep(\al, \be) = (-1)^{(\al, w\be)}\la \al, \be\ra$, but the pairing $\varep_w$ in Section \ref{section-epw}, which is given by $(1 - \zeta)^{(\al, w\be)}(1 - \zeta^2)^{(w\al, \be)}$ in the case $d = 3$, could be used in its place.

\subsection{Construction of a Lie algebra}\label{section-hconstruction}

Assume that $(\Lambda, w, \Heis, \varep)$ is an input datum as defined in the previous section.
Let $T$ be the split torus over $k$ with character group $X^*(T) = \Hom_k(T, \mathbb{G}_m)$ equal to $\Lambda$, and let $\frakt = \text{Lie } T$. We note that $\frakt$ may be identified with $\Hom(\Lambda, k)$, and under this identification $\frakt$ is generated by $\{\al^\vee \mid \al \in \Phi\}$. Let $\mathfrak{l}$ be the $k$-vector space generated by the set of symbols $X_{\wt\al}$ for $\wt\al \in \wt\Phi$, modulo the relation $X_{\zeta\wt\al} = \zeta X_{\wt\al}$ for all $\wt\al \in \wt\Phi$. Thus $\mathfrak{l}$ has dimension $\lvert \Phi \rvert$ as a vector space over $k$. Let $\mathfrak{h} = \mathfrak{t} \oplus \mathfrak{l}$. We define a $k$-bilinear map $[\cdot, \cdot ] : \mathfrak{h} \times \frakh \to \frakh$ by setting
\begin{equation*}
[\al^\vee, X_{\wt\be}] = -[X_{\wt\be}, \al^\vee] = (\al, \be)X_{\wt\be}
\end{equation*}
and
\begin{eqnarray*}
[X_{\wt\al}, X_{\wt\be}] &=&
\begin{cases}
\varep(\al, \be)\wt{\al}\wt{\be}\al^\vee & \text{ if } \al + \be = 0\\
\varep(\al, \be)X_{\wt{\al}\wt{\be}} & \text{ if } \al + \be \in \Phi\\
0 & \text{ otherwise}.
\end{cases}
\end{eqnarray*}

It's not hard to check that this map is well defined. 

\begin{theorem}\label{thm-liealg}
With bracket $[\cdot, \cdot]$ as defined above, the vector space $\mathfrak{h}$ forms a semisimple Lie algebra over $k$ with Cartan subalgebra $\frakt$ and corresponding root lattice $\Lambda$. 
\end{theorem}

\textit{Proof}. First note that $[\cdot, \cdot]$ is alternating. Indeed, if $\al + \be = 0$, then
\begin{eqnarray*}
[X_{\wt\be}, X_{\wt\al}] = \varep(\be, \al)\wt{\be}\wt{\al}\be^\vee
= -\varep(\be, \al)\wt\be\wt\al\al^\vee
= -{[}X_{\wt\al}, X_{\wt\be}{]}
\end{eqnarray*}
since $\wt\al\wt\be = \wt\be\wt\al$ and $\varep(\be, \al) = \varep(-\al, \al) = \varep(\al, \be)$. 
And if $\al + \be \in \Phi$, then 
\begin{equation*}
[X_{\wt\be}, X_{\wt\al}] = \varep(\be, \al)X_{\wt\be\wt\al} = \varep(\be, \al)\la\be, \al\ra X_{\wt\al\wt\be} = -{[}X_{\wt\al}, X_{\wt\be}{]}
\end{equation*}
by the first property of $\varep$. 

We next show that $[\cdot, \cdot]$ satisfies the Jacobi identity. We thus consider
\begin{equation}\label{eqn-bracket}
[x, [y, z]] + [y, [z, x]] + [z, [x, y]]
\end{equation}
for generators $x, y, z$ of $\mathfrak{h}$. If any two of $x, y, z$ are elements of $\mathfrak{t}$, then it's easy to see that (\ref{eqn-bracket}) is zero. Now suppose $x = \al^\vee, y = X_{\wt\be}, z = X_{\wt\gamma}$ for some $\al \in \Phi$ and $\wt\be, \wt\gamma \in \wt\Phi$. If $(\be, \gamma) \geq 0$, then all three terms of (\ref{eqn-bracket}) are zero. If $(\be, \gamma) < 0$, then
(\ref{eqn-bracket}) becomes
\begin{equation*}
(\al, \be + \gamma)[y, z] - (\al, \gamma)[y, z] - (\al, \be)[y, z],
\end{equation*}
which is also zero. Similarly, whenever one of $x, y, z$ is in $\mathfrak{t}$, (\ref{eqn-bracket}) is zero. 

Thus   
from now on we may assume that $x = X_{\wt\al}, y = X_{\wt\be}$, and $z = X_{\wt\gamma}$ for some $\wt\al, \wt\be, \wt\gamma \in \wt\Phi$. 
Suppose $\al + \be + \gamma = 0$.
Then (\ref{eqn-bracket}) is equal to
\begin{eqnarray*}
& &\varep(\be, \gamma)[X_{\wt\al}, X_{\wt\be\wt\gamma}] + \varep(\gamma, \al)[X_{\wt\be}, X_{\wt\gamma\wt\al}] + \varep(\al, \be)[X_{\wt\gamma}, X_{\wt\al\wt\be}]\\
&=& \varep(\be, \gamma)\varep(\al, \be + \gamma)
\wt\al\wt\be\wt\gamma\al^\vee 
+ \varep(\gamma, \al)\varep(\be, \gamma + \al)\wt\be\wt\gamma\wt\al\be^\vee
+ \varep(\al,\be)\varep(\gamma, \al + \be)\wt\gamma\wt\al\wt\be\gamma^\vee\\
&=& \varep(\be, -\al - \be)\varep(\al, -\al)\wt\al\wt\be\wt\gamma\al^\vee
+ \varep(-\al - \be, \al)\varep(\be, -\be)\wt\al\wt\be\wt\gamma\be^\vee\\ 
& &+ \varep(\al, \be)\varep(-\al - \be, \al + \be)\wt\al\wt\be\wt\gamma(-\al^\vee - \be^\vee)\\
&=& \wt\al\wt\be\wt\gamma
[(\varep(\be, -\al - \be)\varep(\al, -\al) - \varep(\al, \be)\varep(-\al -\be, \al + \be))\al^\vee\\
& &+ (\varep(-\al - \be, \al)\varep(\be, -\be) - \varep(\al, \be)\varep(-\al - \be, \al + \be))\be^\vee].
\end{eqnarray*}
Note that
\begin{eqnarray*}
& &\varep(\be, -\al - \be)\varep(\al, -\al) - \varep(\al, \be)\varep(-\al - \be, \al + \be)\\
&=& \varep(\be, -\al)\varep(\be, -\be)\varep(\al, -\al) - \varep(\al, \be)\varep(-\al, \al)\varep(-\al, \be)\varep(-\be, \al)\varep(-\be, \be)\\
&=& \varep(\be, -\be)\varep(\al, -\al)\varep(\be, -\al)(1 - \varep(\al, \be)\varep(-\al, \be))\\
&=& 0
\end{eqnarray*}
by bilinearity of $\varep$, making (\ref{eqn-bracket}) zero.

Thus we can asssume $\al + \be + \gamma \neq 0$. For (\ref{eqn-bracket}) to be nonzero, then at least one of the three terms must be nonzero. Without loss of generality we may assume $[x, [y, z]] = [X_{\wt\al}, [X_{\wt\be}, X_{\wt\gamma}]] \neq 0$, so in particular $[X_{\wt\be}, X_{\wt\gamma}] \neq 0$ and $(\be, \gamma) \in \{-2, -1\}$. We break into two cases depending on the value of $(\be, \gamma)$. 

Case 1: $(\be, \gamma) = -2$, i.e. $\be + \gamma = 0$.

In this case the first term of (\ref{eqn-bracket}) is $-\varep(\be, \gamma)\wt\be\wt\gamma(\be, \al)X_{\wt\al}$.
If $(\al, \be) = -1$, then (\ref{eqn-bracket}) is equal to
\begin{eqnarray*}
\varep(\be, \gamma)\wt\be\wt\gamma X_{\wt\al} + \varep(\al, \be)\varep(\gamma, \al + \be)X_{\wt\gamma\wt\al\wt\be} &=&
(\varep(\be, \gamma)\wt\be\wt\gamma + \varep(\al, \be)\varep(\gamma, \al + \be)\la\al, \be\ra\wt\gamma\wt\be) X_{\wt\al}\\
&=& \varep(\be, -\be)(\wt\be\wt\gamma + \varep(\al, \be)\varep(-\be, \al)\la\al, \be\ra\wt\gamma\wt\be)X_{\wt\al},
\end{eqnarray*}
which is zero by the first property of $\varep$. 

If $(\al, \be) = -2$, then $\al = \gamma = -\be$, so (\ref{eqn-bracket}) is equal to
\begin{eqnarray*}
2\varep(\be, \gamma)\wt\be\wt\gamma X_{\wt\al} + \varep(\al, \be)\wt\al\wt\be[X_{\wt\gamma},\al^\vee] 
&=& 2\varep(\be, \gamma)\wt\be\wt\gamma X_{\wt\al} - 2\varep(\al, \be)\wt\al\wt\be X_{\wt\gamma}\\
&=& 2\varep(\be, \al)\wt\be\wt\gamma X_{\wt\al} - 2\varep(\al, \be)\wt\al\wt\be X_{\wt\gamma\wt\al^{-1}\wt\al}\\
&=& 2\varep(\be, \al)\wt\be\wt\gamma X_{\wt\al} - 2\varep(\al, \be)\wt\al\wt\be\wt\gamma\wt\al^{-1} X_{\wt\al},
\end{eqnarray*}
where we are using the fact that $\wt\gamma\wt\al^{-1} \in \la \zeta\ra$. 
Since $\wt\al$ commutes with $\wt\be$ and $\wt\gamma$ commutes with $\wt\al^{-1}$, we have that $\wt\al\wt\be\wt\gamma\wt\al^{-1} = \wt\be\wt\gamma$, and (\ref{eqn-bracket}) is again zero.

If $(\al, \be) \in \{1, 2\}$, then $(\al, \gamma) = -(\al, \be)$, and (\ref{eqn-bracket}) is equal to $-[X_{\wt\al}, [X_{\wt\gamma}, X_{\wt\be}]] - [X_{\wt\be}, [X_{\wt\al}, X_{\wt\gamma}]]$, which is equal to zero by the cases just considered.

Case 2: $(\be, \gamma) = -1$, i.e. $\be + \gamma \in \Phi$.

In this case, we have $[X_{\wt\al}, [X_{\wt\be}, X_{\wt\gamma}]] = \varep(\be, \gamma)[X_{\wt\al}, X_{\wt\be}]$. Since $\al + \be + \gamma \neq 0$, for the first term (\ref{eqn-bracket}) to be nonzero we must have $(\al, \be + \gamma) = -1$. Thus the first term of (\ref{eqn-bracket}) is $\varep(\al, \be + \gamma)\varep(\be, \gamma)X_{\wt\al\wt\be\wt\gamma} = \varep(\al, \be)\varep(\al, \gamma)\varep(\be, \gamma)X_{\wt\al\wt\be\wt\gamma}$. 

Suppose $(\al, \be)= -1$. Then $(\al, \gamma) = 0$ and $(\al + \be, \gamma) = -1$, and so (\ref{eqn-bracket}) is equal to
\begin{eqnarray*}
\varep(\al, \be)(\varep(\be, \gamma)\varep(\al, \gamma)X_{\wt\al\wt\be\wt\gamma} + \varep(\gamma, \al)\varep(\gamma, \be)X_{\wt\gamma\wt\al\wt\be})
&=& \varep(\al, \be)(\varep(\be, \gamma)\varep(\al, \gamma)+ \varep(\gamma, \al)\varep(\gamma, \be)\la \gamma, \al + \be\ra)X_{\wt\al\wt\be\gamma}\\ &=& 0
\end{eqnarray*}
by the first property of the pairing $\varep$. If $(\al, \gamma) = -1$, the proof is very similar. Otherwise we must have $(\al, \be) = -2$ or $(\al, \gamma) = -2$, and the proof reduces to that of Case 1.

Thus $\frakh$ forms a Lie algebra. It is not hard to show that the Killing form on $\frakh$ is nondegenerate
and that $\frakt$ is a maximal toral subalgebra. The result follows.
\qed

Note that in particular because the Cartan subalgebra $\frakt$ is split, the isomorphism class of $\frakh$ depends only on $\Lambda$ and not on $w, \Heis$, or $\varep$. 

\subsection{Functorality}\label{section-functor}

The construction of $\frakh$ is functorial in the following sense. We define an isomorphism $(\Lambda, w, \Heis, \varep) \to (\Lambda', w', \Heis', \varep')$ of input data to be a pair $(\psi, \phi)$ where 
\begin{itemize}
\item[1.] $\psi: \Lambda \to \Lambda'$ is an isomorphism of simply laced root lattices such that $\psi w = w' \psi$ and $\varep'(\psi(\al), \psi(\be)) = \varep(\al, \be)$ for all $\al, \be \in \Lambda$.  
\item[2.] $\phi: \Heis \to \Heis'$ is an isomorphism making the following diagram commute
\[\begin{tikzcd}        
    1 \arrow{r} & \la\zeta\ra \arrow{r}\arrow[d, "="] & \Heis \arrow{r}\arrow[d, "\phi"] & \Lambda_w\arrow{r}\arrow[d, "\overline\psi"] & 1\\   
  1 \arrow{r} & \la\zeta\ra \arrow{r} & \Heis' \arrow{r} & \Lambda'_{w'} \arrow{r} & 1,
\end{tikzcd}\]
where $\overline\psi$ is the map naturally induced by $\psi$. 
\end{itemize}

Suppose $(\psi, \phi)$ is an isomorphism $(\Lambda, w, \Heis, \varep) \to (\Lambda', w', \Heis', \varep')$ of input data. Let $\frakh$ and $\frakh'$ be the Lie algebras obtained by using the input data $(\Lambda, w, \Heis, \varep)$ and $(\Lambda', w', \Heis', \varep')$, respectively, in the construction of Section \ref{section-hconstruction}. Then $(\psi, \phi)$ naturally induces an isomorphism $\wt\phi: \frakh \to \frakh'$ as follows: with the notation of Section \ref{section-hconstruction}, we have that $\frakh$ is generated by $\{\al^\vee \mid \al \in \Phi\}$ and $\{X_{\wt\al} \mid \wt\al \in \wt\Phi\}$, and $\frakh'$ is generated by $\{\al^\vee \mid \al \in \Phi'\}$ and $\{X'_{\wt\al} \mid \wt\al \in \Phi'\}$, where $\Phi$ and $\Phi'$ are the roots of $\Lambda$ and $\Lambda'$ respectively. 
Note that because $\psi$ and $\phi$ induce the same map on $\Lambda_w'$, they give a well-defined isomorphism (which we will also denote $\phi$) $\wt\Lambda \to \wt\Lambda'$ whose restriction to $\Lambda$ is $\psi$ and whose restriction to $\Heis$ is $\phi$. 
\begin{Lem}\label{lem-functor}
Define $\wt\phi: \frakh \to \frakh'$ by setting
\begin{eqnarray*}
\wt\phi(\al^\vee) &=& \psi(\al)^\vee\\
\wt\phi(X_{\wt\al}) &=& X'_{\phi(\wt\al)}
\end{eqnarray*}
and extending linearly.
Then $\wt\phi$ is a well-defined isomorphism of Lie algebras.
\end{Lem}

Proof.
Note that $\wt\phi$ is well defined because $\phi: \wt\Lambda \to \wt\Lambda$ restricts to the identity on $\la \zeta\ra$. 
It's clear that $\wt\phi$ is a vector-space isomorphism $\frakh \to \frakh'$. It thus suffices to show that $[\wt\phi(x), \wt\phi(y)] = \wt\phi([x,y])$ for all $x, y \in \frakh$. If $x, y \in \frakt$, then it is clear that $\wt\phi([x,y]) = [\wt\phi(x)\wt\phi(y)] = 0$. 
If $\al \in \Phi, \wt\be \in \wt\Phi$, we have that
\begin{eqnarray*}
[\wt\phi(\al^\vee), \wt\phi(X_{\wt\be})] &=& [\psi(\al)^\vee, X'_{\phi(\wt\be)}]\\
&=& (\psi(\al), \psi(\be))X'_{\phi(\wt\be)}\\
&=& (\al, \be)X'_{\phi(\wt\be)}\\
&=& \wt\phi([\al^\vee, X_{\wt\be}]).
\end{eqnarray*}
Now suppose $\wt\al, \wt\be \in \wt\Phi$. If $\al + \be = 0$, we have
\begin{eqnarray*}
[\wt\phi(X_{\wt\al}), \wt\phi(X_{\wt\be})] &=& \varep'(\psi(\al), \psi(\be))\phi(\wt\al)\phi(\wt\be)\psi(\al)^\vee\\
&=& \varep(\al, \be)\phi(\wt\al\wt\be)\psi(\al)^\vee\\
&=& \varep(\al, \be)\wt\al\wt\be\psi(\al)^\vee\\
&=& \wt\phi([X_{\wt\al}, X_{\wt\be}])
\end{eqnarray*}
(we are again using the fact that $\phi$ is the identity on $\wt\al\wt\be \in \la \zeta \ra$). If $(\al, \be) = -1$, then
\begin{eqnarray*}
[\wt\phi(X_{\wt\al}), \wt\phi(X_{\wt\be})] &=& \varep'(\psi(\al), \psi(\be))X'_{\phi(\wt\al\wt\be)}\\
&=& \varep(\al, \be)\wt\phi(X_{\wt\al\wt\be})\\
&=& \wt\phi([X_{\wt\al}, X_{\wt\be}]).
\end{eqnarray*}
The result follows.\qed

\subsection{Defining the grading}\label{section-lieaut}

Now we show that the automorphism $w$ of $\Lambda$ lifts to an order-$d$ automorphism of $\frakh$. 
Note that $w$ extends naturally to an automorphism of
$\wt\Lambda = \Lambda \times_{\Lambda_w} \Heis$ fixing $\Heis$. 
We will also denote this automorphism by $w$ without risk of confusion. 
Define $\lieaut: \frakh \to \frakh$ by setting
\begin{eqnarray*}
\lieaut(\al^\vee) &=& {w(\al)^\vee} \text{ for all } \al \in \Phi,\\
\lieaut (X_{\wt\al}) &=& X_{w(\wt\al)} \text{ for all } \wt\al\in \wt\Phi
\end{eqnarray*} 
and extending linearly.

\begin{Prop}\label{prop-extendw}
The map $\lieaut$ is a Lie algebra automorphism.
\end{Prop}
\textit{Proof}. Let $i$ be the identity automorphism of $\Heis$. Then the pair $(w, i)$ defines an automorphism of input data, as defined in Section \ref{section-functor}. Thus the proposition follows from Lemma \ref{lem-functor}. \qed

It is clear from the definition that $\lieaut$ has order $d$. Let $\frakh_j$ be the $\zeta^j$-eigenspace of $\lieaut$. Thus we have produced a graded Lie algebra 
\begin{equation*}
\frakh = \bigoplus_{j \in \ZZ/d\ZZ} \frakh_j
\end{equation*}
as claimed.

For future considerations, we give a more explicit description of the fixed-point subalgebra $\mathfrak{g} := \frakh^{\lieaut}$. Given $\wt\al \in \wt\Phi$, let
\begin{equation*}
Z_{\wt\al} = X_{\wt\al} + X_{w\wt\al} + ... + X_{w^{d-1}\wt\al}.
\end{equation*}
Then $\mathfrak{g}$ is generated by $\{Z_{\wt\al} \mid \wt\al \in \wt\Phi\}$ (we are using the fact that $w$ is elliptic, which implies that the restriction of $\lieaut$ to $\frakt$ contains no nonzero fixed vectors). 
Note that $Z_{\wt\al} = Z_{w^j\wt\al}$ for all $j$ and that, because $w$ is the identity on $\la \zeta\ra$, we have $Z_{\zeta\wt\al} = \zeta Z_{\wt\al}$ for all $\wt\al\in \wt\Phi$.

We next note that the grading behaves well under functorality. Let $(\psi, \phi)$ be an isomorphism of input data $(\Lambda, w, \Heis, \varep) \to (\Lambda', w', \Heis', \varep')$, and let $\wt\phi: \frakh \to \frakh'$ be the corresponding isomorphism of Lie algebras as defined in Lemma \ref{lem-functor}. Let $\frakg = \frakh^{\lieaut}$ and $\frakg' = \frakh^{\lieaut'}$.

\begin{Prop}\label{prop-functorg}
With notation as above, $\wt\phi$ restricts to an isomorphism $\frakg \to \frakg'$ of Lie algebras. Further, for all $j$, $\wt\phi$ restricts to an isomorphism of vector spaces $\frakh_j \to \frakh_j'$ intertwining the action of $\frakg$ with the action of $\frakg'$. 
\end{Prop}

\textit{Proof}. 
We first claim that $\wt\phi\wt w = \wt w'\wt\phi$. Indeed, if $\al \in \Phi$, we have that
\begin{eqnarray*}
\wt\phi(\wt w(\al^\vee)) &=& \psi (w(\al))^\vee\\
&=& w'(\psi(\al)^\vee)\\
&=& \wt w' \wt\phi (\al^\vee).
\end{eqnarray*}
And for $\wt\al \in\wt\Phi$, we have $\phi(w(\wt\al)) = w'\phi(\wt\al)$, and thus $\wt\phi(\wt w(X_{\wt\al})) = \wt w'(\wt\phi(X_{\wt\al}))$, which proves the claim. The fact that $\wt\phi$ restricts to an isomorphism $\frakh_j \to \frakh'_j$ for all $j$ follows easily. 
In particular, $\wt\phi$ restricts to an isomorphism of subalgebras $\frakg = \frakh_0 \to \frakh'_0 = \frakg'$.
Because $\frakg$ and $\frakg'$ act via restriction of the adjoint representations of $\frakh$ and $\frakh'$ respectively, the fact that $\wt\phi$ intertwines the actions of $\frakg$ and $\frakg'$ is equivalent to the fact that $\wt\phi$ is a Lie algebra isomorphism.
\qed

\begin{remark}\label{rem-innerauts}
Applying Lemma \ref{lem-functor} to the inner automorphisms of $\Heis$, we obtain a natural homomorphism
\begin{equation*}
\iota: \Lambda_w \to \Aut(\frakh).
\end{equation*}
We have that $\iota(\lam)(X_{\wt\al}) = \la\lam, \al\ra X_{\wt\al}$ for all $\lam \in \Lambda_w, \wt\al \in \wt\Phi$, so $\iota$ is injective whenever $\la \cdot, \cdot \ra$ is nondegenerate. 
By Proposition \ref{prop-functorg}, the elements of $\iota(\Lambda_w)$ restrict to automorphisms of $\frakg$, yielding a homomorphism 
\begin{equation*}
\Lambda_w \to \Aut(\frakg).
\end{equation*}
\end{remark}

\subsection{Galois descent}\label{section-galois}

For number-theoretic considerations, we are often interested in a more general version of the construction of the previous sections, in which 
we replace extensions of the form (\ref{eqn-centralext}) with certain central extensions of \'etale group schemes (see \cite{Tho16}, \cite{RomanoThorne2}). 
In particular, we may consider a $\mu_d$-extension of $J_C[d]$,
where $J_C$ is the Jacobian of a curve $C$ defined over  
some number field not necessarily containing a primitive $d$th root of unity (Section \ref{section-examples} below discusses central extensions of Jacobians further). Because of this, it becomes important to show that we may use the construction of the previous sections to obtain Lie algebras over arbitrary fields of characteristic zero.
 The next proposition shows that under certain conditions, this becomes possible with Galois descent.

Let $k$ be as above, let $(\Lambda, w, \Heis, \varep)$ be an input datum, and let $\frakh$ and $\frakg = \frakh^{\wt w}$ be the associated Lie algebras as constructed above. Suppose $k/k_0$ is a finite Galois extension, and let $\Gamma = \Gal(k/k_0)$. Define a homomorphism $\theta: \la \zeta\ra \to \Aut(\Lambda)$ by $\theta(\zeta) = w$. 
Suppose we are given homomorphisms $\Gamma \to \Aut(\Lambda)$ and $\Gamma \to \Aut(\Heis)$ such that 
\begin{enumerate}
\item[A1.] $\sigma\theta(\zeta^i)\sigma^{-1} = \theta(\sigma(\zeta^i))$ in $\Aut(\Lambda)$ for all $\sigma \in \Gamma$ and all $i$.
\item[A2.] The action of $\Gamma$ on $\Heis$ restricts to the natural Galois action on $\la \zeta \ra$. 
\item[A3.] The naturally induced homomorphisms $\Gamma \to \Aut(\Lambda) \to \Aut(\Lambda_w)$ and $\Gamma \to \Aut(\Heis) \to \Aut(\Lambda_w)$  are the same.
\end{enumerate}

In this setting, the following proposition gives sufficient conditions to define forms of $\frakh$ and $\frakg$ over $k_0$. We note that the two conditions of Proposition \ref{prop-basefield} are automatically satisfied in the case when $\la \cdot, \cdot \ra = \la \cdot, \cdot \ra_w$ and $\varep = \varep_w$ as defined in Section \ref{section-pairingsandexts} (see Lemma \ref{lem-epwgalois}).

\begin{Prop}\label{prop-basefield}
With notation as above, suppose the following two conditions hold for all $\sigma \in \Gamma$:
\begin{enumerate}
\item[1.] $\sigma (\la \al, \be \ra) = \la \sigma(\al), \sigma(\be)\ra$ for all $\al, \be \in \Lambda$.
\item[2.] $\sigma(\varep(\al, \be)) = \varep(\sigma(\al), \sigma(\be))$ for all $\al, \be \in \Lambda$.
\end{enumerate}
Then there is a natural semilinear action of $\Gamma$ on $\frakh$ preserving $\frakg$ such that $\frakh^\Gamma$ and $\frakg^\Gamma$ are 
$k_0$-forms of the Lie algebras $\frakh$ and $\frakg$, respectively.
\end{Prop}

\textit{Proof}. 
Because the $\Gamma$-actions on $\Lambda$ and $\Heis$ induce the same action on $\Lambda_w$, these actions naturally induce a $\Gamma$-action on $\wt\Lambda$.
Given $\sigma \in \Gamma$, define a $\sigma$-linear map $a_\sigma: \frakh \to \frakh$ by $a_\sigma(\al^\vee) = \sigma(\al)^\vee$ and $a_\sigma (X_{\wt\al}) = X_{\sigma(\wt\al)}$ for all $\wt\al \in \wt\Phi$. Note $a_\sigma$ is well defined because the action of $\Gamma$ on $\wt\Lambda$ restricts to the natural Galois action on $\la \zeta \ra$. 
We have that
\begin{eqnarray*}
[a_\sigma(\al^\vee), a_\sigma(\be^\vee)] = a_\sigma([\al^\vee, \be^\vee]) = 0 
\end{eqnarray*}
for all $\al, \be \in \Lambda$;
\begin{eqnarray*}
[a_\sigma(\al^\vee), a_\sigma(X_{\wt\be})] &=& (\sigma(\al), \sigma(\be))X_{\sigma(\wt\be)}\\
&=& (\al, \be)X_{\sigma(\wt\be)}\\
&=& a_\sigma([\al^\vee, X_{\wt\be}])
\end{eqnarray*}
for all $\al \in \Lambda, \wt\be \in \wt\Phi$;
\begin{eqnarray*}
[a_\sigma(X_{\wt\al}), a_\sigma(X_{\wt\be})] &=& \varep(\sigma(\al), \sigma(\be))\sigma(\wt\al)\sigma(\wt\be)\sigma(\al)^\vee\\
&=& \sigma(\varep(\al, \be))\sigma(\wt\al\wt\be)\sigma(\al)^\vee\\
&=& a_\sigma([X_{\wt\al}, X_{\wt\be}])
\end{eqnarray*}
if $\wt\al, \wt\be \in \wt\Phi$ with $\al + \be = 0$; and
\begin{eqnarray*}
[a_\sigma(X_{\wt\al}), a_\sigma(X_{\wt\be})] &=& \varep(\sigma(\al), \sigma(\be))X_{\sigma(\wt\al)\sigma(\wt\be)}\\
&=& a_\sigma([X_{\wt\al}, X_{\wt\be}]
\end{eqnarray*}
if $\wt\al, \wt\be \in \wt\Phi$ with $(\al, \be) = -1$.
It follows 
that $a_\sigma$ is a $k_0$-algebra automorphism, and the collection of maps $\{a_\sigma \mid \sigma \in \Gamma\}$ defines a $\Gamma$-structure on $\frakh$.  Thus $\frakh^\Gamma$ is a $k_0$-form of $\frakh$.
Also note that $a_\sigma(Z_{\wt\al}) = Z_{\sigma(\wt\al)}$ for all $\wt\al \in \wt\Phi$, so the action of $\Gamma$ preserves $\frakg$, and thus $\frakg^\Gamma$ is a $k_0$-form of $\frakg$. \qed

Note that if we define stacks in the same way as \cite[Section 3]{RomanoThorne2}, then the proofs above would be enough to define a morphism of stacks as in \cite[Theorem-Construction 3.1]{RomanoThorne2}.

\section{Bilinear pairings and central extensions}\label{section-pairingsandexts}

Let $\Lambda$ be a simply laced root lattice, and let $w \in \Aut(\Lambda)$ be elliptic.
In this section we show that given $\Lambda$ and $w$, we may explicitly realize an input datum as defined in Section \ref{section-input}.

\subsection{The pairing $\la \cdot, \cdot \ra_w$}\label{section-pairing}

We first describe a bilinear pairing on $\Lambda_w$. We begin with two pairings on $\Lambda$: the first, $C$, was defined by Lepowsky in \cite[Section 4]{Lepowsky}; the second, $\la \cdot, \cdot \ra_w$ was defined by Reeder in \cite{ReederEll}\footnote{Note that \cite{ReederEll} gives an additive version of this pairing, but we require the associated multiplicative pairing described here.} and was shown to factor through $\Lambda_w \times \Lambda_w$. We will then show that these two pairings are the same.

The pairing $C: \Lambda \times \Lambda \to k^\times$ of \cite{Lepowsky} is defined by
\begin{equation*}
C(\al, \be)= \zeta^{(\sum_{j = 1}^{d - 1} jw^j\al, \be)}.
\end{equation*} 
Now let $M(t)$ be the minimal polynomial of $w$. By \cite[Lemma 2.1]{ReederEll}, $m:= M(1)$ divides $d$; we will write $d_0$ for $\frac{d}{m}$. Let $\zeta_m = \zeta^{d_0}$, which is a primitive $m$th root of unity. 
Let
\begin{equation*}
{M_0}(t) = \frac{M(t) - M(1)}{t - 1},
\end{equation*}
which is a polynomial over $\mathbb{Z}$. Given $\al, \be \in \Lambda$, we set
\begin{equation*}
\langle \al, \be\rangle_w = \zeta_m^{( \al, {M_0}(w)\be)}.
\end{equation*}
Note that $\la \al, \be \ra_w$ depends only on the images of $\al$ and $\be$, respectively, in $\Lambda_w$, and the pairing $\la \cdot, \cdot \ra_w$ is alternating (\cite[Lemma 2.2]{ReederEll}).

\begin{Lem}\label{lem-2pairings} We have that
\begin{equation*}
\la\al, \be \ra_w = C(\al, \be)
\end{equation*} 
for all $\al, \be \in \Lambda$.
\end{Lem}

\textit{Proof}. 
Let $g(t) = 1 + t + t^2 + ... + t^{d - 1}$, so that $C(\al, \be) = \zeta^{(wg'(w)\al, \be)}$. 
We have $g(t)(t - 1) = t^d - 1$, so
\begin{equation*}
g'(t)(t - 1) = dt^{d - 1} - g(t), 
\end{equation*}
and
\begin{equation*}
tg'(t)(t - 1) = dt^d - tg(t).
\end{equation*}
 Plugging in $w$, we have
$wg'(w)(w - 1) = d$. Since $w - 1$ is an invertible linear transformation of $\Lambda \otimes_{\ZZ} \mathbb{R}$, the inverse $(w - 1)^{-1}$ is well defined, and we have 
\begin{equation*}
wg'(w) = d(w - 1)^{-1}.
\end{equation*}
Similarly, we have
$M_0(w)(w - 1) = -m$, and $M_0(w) = -m(w - 1)^{-1}$. Thus we see that 
\begin{equation*}
wg'(w) = -d_0M_0(w).
\end{equation*} 
Hence
\begin{eqnarray*}
C(\al, \be) &=& \zeta^{(wg'(w)\al, \be)}\\
&=& \zeta_m^{-(M_0(w)\al, \be)}\\
&=& \la \be, \al\ra_w^{-1}\\
&=& \la \al, \be\ra_w.
\end{eqnarray*}
\qed

As a consequence, the pairing $C$ factors through $\Lambda_w \times \Lambda_w$. We will 
now drop the notation $C$ and just write $\la \cdot, \cdot \ra_w$ for the pairing of Lemma \ref{lem-2pairings}.

\subsection{The existence of central extensions}

We now review some of the theory of central extensions to show that there exists a central extension of $\Lambda_w$ by $\la \zeta \ra$ with commutator pairing $\la\cdot, \cdot \ra_w$. For a nice overview of this topic see \cite{Prasad}, and for more details see, e.g., \cite[Chapter 1]{BeylTappe}.

Isomorphism classes of central extensions
\begin{equation*}
1 \to \la \zeta \ra \to E \to \Lambda_w \to 1
\end{equation*}
are parametrized by $H^2(\Lambda_w, \la \zeta \ra)$, where $\la \zeta \ra$ is considered to be a trivial $\Lambda_w$-module. Recall that there is a short exact sequence
\begin{equation*}
 0 \to \text{Ext}^1(\Lambda_w, \la \zeta \ra) \to H^2(\Lambda_w, \la \zeta \ra) \to \Hom(\wedge^2 \Lambda_w, \la \zeta \ra) \to 0,
\end{equation*}
where we are taking $\text{Ext}$ in the category of abelian groups, and $\wedge^2 \Lambda_w = (\Lambda_w \otimes \Lambda_w)/\{x \otimes x \mid x \in \Lambda_w\}$.
We may view the map $H^2(\Lambda_w, \la \zeta \ra) \to \Hom(\wedge^2 \Lambda_w, \la \zeta \ra)$ as sending a central extension to its commutator pairing. 
Thus because the bilinear pairing $\la \cdot, \cdot \ra_w : \Lambda_w \times \Lambda_w \to \la \zeta \ra$ is alternating, there exists a central extension of $\Lambda_w$ by $\la \zeta \ra$ with commutator pairing given by $\la \cdot, \cdot \ra_w$. 

More explicitly, a cocycle $c$ defines multiplication on the set $\la \zeta \ra \times \Lambda_w$ by $(\zeta^i, \al)\cdot (\zeta^j, \be) = (\zeta^{i + j}c(\al, \be), \al + \be)$. With this definition, the inverse to $(\zeta^i, \al)$ is $(\zeta^{-1}c(-\al, \al)^{-1}c(0, 0)^{-1}, -\al)$. Using these formulas, 
it's not hard to check that if a cocycle $c$ is bilinear, the commutator pairing corresponding to $c$ is given by $c(\al, \be)c(\be, \al)^{-1}$. 

As an example of a choice of $c$, write $\Lambda_w$ as a direct produce of cyclic groups $C_1 \times ... \times C_n$, and choose a generator $x_i$ for each $C_i$. 
Define a cocyle on $\Lambda_w$ by setting 
\begin{eqnarray*}
c(x_i, x_j) &=&
\begin{cases}
\la x_i, x_j \ra_w \text{ if } i < j\\
1 \text{ if } i \geq j
\end{cases}
\end{eqnarray*}
and extending bilinearly. Then, using the fact that $\la \cdot, \cdot \ra_w$ is skew-symmetric, it is not hard to check that $c(\al, \be)c(\be, \al)^{-1} = \la \al, \be \ra_w$ for all $\al, \be \in \Lambda_w$. 

\begin{remark}
Note that if $(\Lambda, w, \Heis, \varep)$ is an input datum, then the commutator pairing on $\Heis$ is not necessarily given by $\la \cdot, \cdot \ra_w$. 
In fact, if $\Lambda$ is a simply laced root lattice and $w \in \Aut(\Lambda)$ has odd order, then given any alternating pairing $\la\cdot, \cdot\ra$ on $\Lambda_w$, there exists an input datum $(\Lambda, w, \Heis, \varep)$ such that $\Heis$ has commutator pairing $\la \cdot, \cdot \ra$. Indeed, 
by the logic above, there exists a central extension $\Heis$ with commutator pairing given by $\la\cdot, \cdot\ra$. If we define $\varep: \Lambda \times \Lambda \to k^\times$ by 
\begin{equation*}
\varep(\al, \be) = (-1)^{(\al, \sum_{j = 1}^{\frac{d-1}{2}} w^j\be)}\la \al, \be\ra^{\frac{d - 1}{2}},
\end{equation*}
then it's not hard to check that $(\Lambda, w, \Heis, \varep)$ is an input datum.
\end{remark}

\subsection{The pairing $\varep_w$}\label{section-epw}

We now define a pairing $\varep_w: \Lambda \times \Lambda \to k^\times$ satisfying the two conditions in the definition of an input datum with respect to $\la \cdot, \cdot \ra_w$. We note that this pairing appears in \cite[Section 4]{Lepowsky}. Let 
\begin{equation*}
\varep_w(\al, \be) = \prod_{j = 1}^{d - 1} (1 - \zeta^{-j})^{(w^j\al, \be)}.
\end{equation*} 
It is clear that $\varep_w(w\al, w\be) = \varep_w(\al, \be)$ for all $\al, \be \in \Lambda$. 
The following lemma, which is stated but not proven in \cite[Section 4]{Lepowsky}, shows that $\varep_w$ satisfies the necessary properties described in Section \ref{section-input}.

\begin{Lem}
We have
\begin{equation*}
\frac{\varep_w(\al, \be)}{\varep_w(\be, \al)} = (-1)^{(\al, \be)}\la\be, \al\ra_w
\end{equation*}
for all $\al, \be \in \Lambda$. 
\end{Lem}

\textit{Proof}. We have
\begin{eqnarray*}
\frac{\varep_w(\al, \be)}{\varep_w(\be, \al)} &=& \frac{\prod_{j = 1}^{d - 1} (1 - \zeta^{-j})^{(w^j\al, \be)}}{\prod_{j = 1}^{d - 1} (1 - \zeta^{-j})^{(w^j\be, \al)}}\\
&=& \frac{\prod_{j = 1}^{d - 1} (1 - \zeta^{-j})^{(w^j\al, \be)}}{\prod_{j = 1}^{d - 1} (1 - \zeta^{j})^{(w^j\al, \be)}}\\
&=& \frac{\prod_{j = 1}^{d - 1} (1 - \zeta^{-j})^{(w^j\al, \be)}}{\prod_{j = 1}^{d - 1} (-\zeta^j)^{(w^j\al, \be)}(1 - \zeta^{-j})^{(w^j\al, \be)}}\\
&=& (-1)^{(\al, \be)} \la \al, \be\ra_w^{-1}.
\end{eqnarray*}
\qed

Thus we see that if $\Heis$ is any central extension of $\Lambda_w$ by $\la \zeta \ra$ with commutator pairing $\la \cdot, \cdot \ra_w$, then $(\Lambda, w, \Heis, \varep_w)$ is an input datum. We next note that an input datum of this form behaves well with respect to functorality and Galois descent (cf. sections \ref{section-functor} and \ref{section-galois}).   

\begin{Lem}\label{lem-epw-functor}
Suppose $\Lambda$ and $\Lambda'$ are simply laced root lattices with elliptic automorphisms $w$ and $w'$ respectively. If $\psi: \Lambda \to \Lambda'$ is an isomorphism of root lattices such that $\psi w = w'\psi$, then 
\begin{equation*}
\varep_{w'}(\psi(\al), \psi(\be)) = \varep_w(\al, \be)
\end{equation*}
for all $\al, \be \in \Lambda$.
\end{Lem}

\textit{Proof}. This follows directly from definitions. \qed

As mentioned earlier, the next result shows that $\la \cdot, \cdot \ra_w$ and $\varep_w$ satisfy the conditions of Proposition \ref{prop-basefield}.

\begin{Lem}\label{lem-epwgalois}
With notation as in Section \ref{section-galois}, suppose $\Gamma$ acts on $\Lambda$ such that $\sigma\theta(\zeta)\sigma^{-1} = \theta(\sigma(\zeta))$ for all $\sigma \in \Gamma$. Then 
\begin{equation*}
\sigma(\la \al, \be \ra_w) = \la \sigma(\al), \sigma(\be) \ra_w
\end{equation*} and 
\begin{equation*}
\sigma(\varep_w(\al, \be)) = \varep_w(\sigma(\al), \sigma(\be))
\end{equation*}
 for all $\al, \be \in \Lambda, \sigma \in \Gamma$. 
\end{Lem}

\textit{Proof}. Suppose $\sigma(\zeta) = \zeta^n$, and let $n^{-1}$ denote the inverse of $n$ in $(\ZZ/d\ZZ)^\times$. Then
\begin{eqnarray*}
\sigma(\la \al, \be \ra_w) &=& (\zeta^n)^{(\sum_{j = 1}^{d-1} jw^j \al, \be)}\\
&=& \zeta^{(\sum_{j = 1}^{d - 1} jnw^j\al, \be)}\\
&=& \zeta^{(\sum_{j = 1}^{d - 1} jw^{jn^{-1}}\al, \be)}\\
&=& \zeta^{(\sum_{j = 1}^{d - 1} jw^{j} \sigma(\al), \sigma(\be)}\\
&=& \la \sigma(\al), \sigma(\be) \ra_w.
\end{eqnarray*}
Here we are using the fact that $\sigma^{-1}w\sigma = w^{n^{-1}}$. The proof is similar for $\varep_w$.\qed

\subsection{On representations of $\Heis$}

Section \ref{section-reps-of-g} will relate the representations of the fixed-point subalgebra $\frakg$ to the representations of $\Heis$, so for completeness, we recall here some facts about representation theory of $\Heis$. In this subsection, $\Lambda, w$, and $\zeta$ are as above, and $\Heis$ may be any central extension of $\Lambda_w$ by $\la \zeta\ra$. 
Let $K$ be an algebraically closed field containing $k$, let $Z$ be the center of $\Heis$, and let $A$ be a maximal abelian subgroup of $\Heis$. The following proposition summarizes some well-known properties.

\begin{Prop}\label{prop-repsofHeis}
Let $\chi: Z \to K^\times$ be a character of $Z$ that restricts to the natural inclusion $\la \zeta \ra \hookrightarrow K^\times$, and let $\wt\chi$ be any extension of $\chi$ to $A$. Then
\begin{itemize}
\item[1.] $\Ind_A^\Heis\wt\chi$ is an irreducible representation on which $\la \zeta \ra$ acts via its defining character.
\item[2.] If $\chi': Z \to K^\times$ is a character restricting to the inclusion $\la \zeta \ra \hookrightarrow K^\times$ and $\wt\chi'$ is an extension of $\chi'$ to $A$, then $\Ind_A^\Heis \wt\chi \simeq \Ind_A^\Heis \wt\chi'$ as representations of $\Heis$ if and only if $\chi = \chi'$.
\item[3.] Every irreducible representation of $\Heis$ over $K$ on which $\la \zeta \ra$ acts by its defining character is isomorphic to one of the form $\Ind_A^\Heis \wt\chi$ for some $\chi$ as above.
\end{itemize}
\end{Prop}

\textit{Proof}. 
The fact that $\la \zeta \ra$ acts via its defining character follows from the definition of an induced representation.
An easy exercise in Mackey theory shows that $\Ind_A^{\Heis} \wt\chi$ is irreducible, and that if $\chi': Z \to k^\times$ is as in part 2 of the proposition, then $\Ind_A^{\Heis}\wt{\chi} \simeq \Ind_A^{\Heis} \wt{\chi}'$ implies $\chi = \chi'$. Now given an irreducible representation $V$ on which $\la \zeta \ra$ acts by its defining character, choose a character $\wt\chi$ of $A$ such that $\Hom_A(\wt\chi, V) \neq 0$. Then by Frobenius reciprocity, we have $\Ind_A^\Heis \wt\chi \simeq V$.\qed

\section{Representations of the fixed-point subalgebra}\label{section-reps-of-g}

Let $(\Lambda, w, \Heis, \varep_w)$ be an input datum containing $\varep_w$ as defined in Section \ref{section-epw}. Let $\frakh$ and $\frakg = \frakh^{\wt w}$ be the corresponding Lie algebras resulting from the construction of  Section \ref{section-construction}, and let $K$ be any field extension of $k$. 
In this section we show that every representation of $\Heis$ over $K$ on which $\la \zeta \ra$ acts by its defining character may be used to define a representation of $\frakg$.

It will be helpful to introduce the following notation. Given $\al, \be \in \Phi$ and $c \in \{0, \pm 1, \pm 2\}$, let
\begin{equation*}
I_{\al, \be}(c) = \{j \in \{0, \dots, d-1\} \mid (w^j\al, \be) = c\}.
\end{equation*}

We begin with a technical result, Proposition \ref{prop-poly-relation}, that will play a key role in the proof of Theorem \ref{thm-rep-ext}. Proposition \ref{prop-poly-relation} follows from \cite[Proposition 4.1]{Lepowsky}, but I have been unable to find a proof of it in the literature. Thus for completeness, I include a proof here.

\begin{Prop}\label{prop-poly-relation}
Suppose $\al, \be \in \Phi$ with $(w^j\al, \be) \in \{0, \pm 1\}$ for all $j$. Then
\begin{equation}\label{eqn-poly-relation}
1 - \la \be, \al \ra_w = \sum_{j \in I_{\al, \be}(-1)} \varep_w(w^j\al, \be).
\end{equation}
\end{Prop}

We begin by stating a lemma, the proof of which is a straightforward exercise.

\begin{Lem}\label{lem-poly-sum}
Let $P(T) = \sum_{i = 0}^n c_iT^i \in K[T]$. 
Then 
\begin{equation*}
\sum_{j = 0}^{d - 1} P(\zeta^j) =  \underset{\substack{i \equiv 0\\ \mod d}
}{\sum} dc_i
\end{equation*}
\end{Lem}

\textit{Proof of Proposition \ref{prop-poly-relation}}.
Fix $\al, \be \in \Phi$ as in the statement of the proposition. We define
\begin{equation*}
P(T) = \prod_{j = 0}^{d - 1} (1 - \zeta^{-j}T)^{(w^j\al, \be) + 1}.
\end{equation*} 
By assumption,
we have $P \in K[T]$. We will show that both sides of (\ref{eqn-poly-relation}) are equal to 
\begin{equation}\label{eqn-poly-sum}
\frac{1}{d}\sum_{j = 0}^{d-1} P(\zeta^j).
\end{equation}
We first calculate the value of (\ref{eqn-poly-sum}) using Lemma \ref{lem-poly-sum}.
Since $(1 + w + w^2 + ... + w^{d-1})\al = 0$, we see that $\lvert I_{\al, \be}(1) \rvert= \lvert I_{\al, \be}(-1) \rvert$, and thus the degree of $P$ is $d$. If we write $P(T) = \sum c_iT^i$, by Lemma \ref{lem-poly-sum} we have that
$\frac{1}{d}\sum_{j = 0}^{d-1} P(\zeta^j) = c_0 + c_d$. It's clear that $c_0 = 1$. We now calculate $c_d$:
\begin{eqnarray*}
c_d &=& \prod_{j = 0}^{d-1} (-\zeta^{-j})^{(w^j\al, \be) + 1}\\
&=& \prod_{j = 0}^{d-1} (-\zeta^{-j})\prod_{j = 0}^{d-1} (-\zeta^{-j})^{(w^j\al, \be)}\\
&=& (-1)^d\prod_{j = 1}^{d - 1} (\zeta^{-j})\la \al, \be\ra_w^{-1}\\
&=& (-1)^d\zeta^{-\sum_{j = 1}^{d-1} j}\la\be,\al\ra_w.
\end{eqnarray*}
We have 
\begin{eqnarray*}
\zeta^{-\sum_{j = 1}^{d-1} j} &=& \zeta^{\frac{d(d-1)}{2}},
\end{eqnarray*}
which is 1 if $d$ is odd and $-1$ if $d$ is even. Thus $c_d = -\la \be, \al \ra_w$, and (\ref{eqn-poly-sum}) is equal to $1 - \la \be, \al \ra_w$, which is the left-hand side of (\ref{eqn-poly-relation}).

We now show that the right-hand side of (\ref{eqn-poly-relation}) is also equal to (\ref{eqn-poly-sum}).
Note that for all $i \in \{0, ..., d- 1\}$, we have $P(\zeta^i) = 0$ if and only if $(w^i\al, \be) \geq 0$. 
If $i \in I_{\al, \be}(-1)$, then 
\begin{eqnarray*}
P(\zeta^i) &=& \prod_{j = 0}^{d-1} (1 - \zeta^{i-j})^{(w^j\al, \be) + 1}\\
&=& \prod_{j \neq i} (1 - \zeta^{i-j})^{(w^j\al, \be) + 1}\\
&=& \prod_{j \neq i}(1 - \zeta^{i-j})^{(w^j\al, \be)}\prod_{j \neq i}(1 - \zeta^{i-j})\\
&=& \prod_{j \neq i}(1 - \zeta^{i-j})^{(w^j\al, \be)}\prod_{j = 1}^{d-1}(1 - \zeta^{j})\\
&=&
d\prod_{j \neq i}(1 - \zeta^{i-j})^{(w^j\al, \be)}\\ 
&=& d\varep_w(w^i\al, \be).
\end{eqnarray*}
Thus the right-hand side of (\ref{eqn-poly-relation}) is also equal to (\ref{eqn-poly-sum}) and we have proved the proposition. \qed

We now consider $\mathfrak{g}$. 
Recall that given $\wt\al \in \wt\Phi$, we defined
\begin{equation*}
Z_{\wt\al} = X_{\wt\al} + X_{w\wt\al} + ... + X_{w^{d-1}\wt\al}.
\end{equation*}

\begin{Lem}\label{lem-Z-bracket}
Given $\wt\al, \wt\be \in \wt\Phi$, we have
\begin{equation*}
[Z_{\wt\al}, Z_{\wt\be}] = \sum_{j \in I_{\al, \be}(-1)} \varep_w(w^j\al, \be)Z_{(w^j\wt\al)\wt\be}.
\end{equation*}
\end{Lem}

Proof. By definition, $[Z_{\tilde\al}, Z_{\tilde\be}]$ is given by
\begin{equation*}
\begin{matrix}
& {[} X_{\wt\al}, X_{\wt\be}{]} & + & {[}X_{w\wt\al}, X_{w\wt\be}{]} & + & {[}X_{w^2\wt\al}, X_{w^2\wt\be}{]} & + & \dots & + & {[}X_{w^{d-1}\wt\al}, X_{w^{d-1}\wt\be}{]}\\
+ & {[}X_{w\wt\al}, X_{\wt\be}{]} &+ & {[}X_{w^2\wt\al}, X_{w\wt\be}{]}& + & \dots & & & + & {[}X_{\wt\al}, X_{w^{d-1}\wt\be}{]}\\
& \vdots & & & & & & & & \\
+ &  {[}X_{w^{d-1}\wt\al}, X_{\wt\be}{]} &+& {[}X_{\wt\al}, X_{w\wt\be}{]} & +&  \dots & & &  + & {[}X_{w^{d-2}\wt\al}, X_{w^{d-1}\wt\be}{]}.
 \end{matrix}
\end{equation*}
If $(w^j\al, \be) \geq 0$, then each term in the sum
\begin{equation*}
[X_{w^j\wt\al}, X_{\wt\be}] + [X_{w^{j + 1}\wt\al}, X_{w\wt\be}] + \dots + [X_{w^{j - 1}\wt\al}, X_{w^{d-1}\wt\be}]
\end{equation*}
is zero. If $(w^j\al, \be) = -2$, then $(w^j\wt\al)\wt\be \in \mu_d$, so $(w^{j + i}\wt\al)(w^i\wt\be) = w^i((w^j\wt\al)\wt\be)) = w^j\wt\al\wt\be$ for all $i$. Thus
\begin{eqnarray*}
[X_{w^j\wt\al}, X_{\wt\be}] &+& [X_{w^{j + 1}\wt\al}, X_{w\wt\be}] + \dots + [X_{w^{j - 1}\wt\al}, X_{w^{d-1}\wt\be}]\\
&=& \varep_w(w^j\al, \be)w^j\wt\al\wt\be(w^j\al^\vee + w^{j + 1}\al^\vee + \dots + \al^\vee + w\al^\vee + \dots + w^{j - 1}\al^\vee)\\
&=& 0.
\end{eqnarray*}
Finally, if $(w^j\al, \be) = -1$, then
\begin{equation*}
[X_{w^j\wt\al}, X_{\wt\be}] + [X_{w^{j + 1}\wt\al}, X_{w\wt\be}] + \dots + [X_{w^{j - 1}\wt\al}, X_{w^{d-1}\wt\be}] = \varep_w(w^j\al, \be)Z_{w^j\wt\al\wt\be}.
\end{equation*}
The result follows. \qed

We now turn to representations of $\Heis$.  Let $(\rho, V)$ be a representation of $\Heis$ over $K$ such that $\rho(\zeta)v = \zeta v$ for all $v \in V$. Define a linear transformation $\tilde\rho: \mathfrak{g} \to \mathfrak{gl}(V)$ as follows: let $\pi: \wt\Lambda = \Lambda \times_{\Lambda_w} \Heis \to \Heis$ be the canonical projection. 
Define $\tilde\rho$ by setting
\begin{equation*}
\tilde\rho(Z_{\wt\al}) = \rho(\pi(\wt\al))
\end{equation*} 
for all $\wt\al \in \wt\Phi$.
Then $\tilde\rho$ is well defined because 
$\pi(w^j\wt\al)= \pi(\wt\al)$ 
for all $j$ and $\wt\al$, and because $\rho(\zeta)$ is scalar multiplication by $\zeta$. 

\begin{theorem}\label{thm-rep-ext}
The map $\tilde\rho: \frakg \to \mathfrak{gl}(V)$ is a Lie algebra homomorphism; i.e. $\tilde\rho$ defines a representation of $\mathfrak{g}$ on $V$.
\end{theorem}

\textit{Proof}. 
It suffices to show that
\begin{equation}\label{eqn-rep-lift}
\tilde\rho([Z_{\wt\al}, Z_{\wt\be}]) = [\tilde\rho(Z_{\wt\al}), \tilde\rho(Z_{\wt\be})]
\end{equation}
for all $\wt\al, \wt\be \in \wt\Phi$. 
We break into three cases.

Case 1: $(\al, w^i\be) = 2$ for some $i$.

Replacing $\be$ by $w^i\be$, we may assume $i = 0$. Then $\al = \be$, so $\wt\be = \zeta^n\wt\al$ for some $n$. Since $Z_{\wt\be} = \zeta^n Z_{\wt\al}$, we see that the left-hand side of (\ref{eqn-rep-lift}) is zero. The right-hand side is 
\begin{eqnarray*}
\rho(\pi(\wt\al)\pi(\wt\be)) - \rho(\pi(\wt\be)\pi(\wt\al)) &=& \rho(\zeta^n\pi(\wt\al)^2) - \rho(\zeta^n\pi(\wt\al)^2)\\
&=& 0.
\end{eqnarray*}
For the rest of the proof we assume that $(\al, w^i\be) \neq 2$ for all $i$.

Case 2: $(\al, w^i\be) = -2$ for some $i$.

By the same logic as above, we may assume $i = 0$. We have that $\be = -\al$ and $\wt\be = \zeta^n\wt\al^{-1}$ for some $n$. Thus $\pi(\wt\al)\pi(\wt\be) \in \mu_d$, making the right-hand side of (\ref{eqn-rep-lift}) zero.
To see that the left-hand side of (\ref{eqn-rep-lift}) is zero, by Lemma \ref{lem-Z-bracket}, it suffices to show that 
\begin{equation*}\label{eqn-I-1}
\sum_{j \in I_{\al, \be}(-1)} \varep_w(w^j\al, \be)\rho(\pi(w^j\wt\al\wt\be))
\end{equation*}
is zero. First note that $\frac{d}{2} \notin I_{\al, \be}(-1)$.
Indeed, 
since $\al + w\al + w^2\al + ... + w^{d-1}\al = 0$, we have
\begin{eqnarray*}
(\al, w\al + w^2\al + ... + w^{d-1}\al) = -2.
\end{eqnarray*}
If $d$ is even, 
this equation becomes
\begin{eqnarray*}
-2 &=& 2(\al, w\al + w^2\al + ... + w^{\frac{d}{2} - 1}\al) + (\al, w^{\frac{d}{2}}\al),
\end{eqnarray*}
which implies that $(\al, w^{\frac{d}{2}}\al) \in 2\mathbb{Z}$. 

Next note that if $(w^j\al, \be) = -1$, then $(w^{d-j}\al, \be) = (\al, w^j\be) = (w^j\al, \be) = -1$ since $\al = -\be$. Thus it suffices to show that if $j \in I_{\al, \be}(-1)$, then 
\begin{equation}\label{eqn-2terms}
\varep_w(w^j\al, \be)\rho(\pi(w^j\wt\al\wt\be)) + \varep_w(w^{d-j}\al, \be)\rho(\pi(w^{d-j}\wt\al\wt\be))
\end{equation}
is zero. Fix $j \in I_{\al, \be}(-1)$. Note that since $\al = -\be$, we have
\begin{equation*}
\varep_w(w^{d-j}\al, \be) 
= \varep_w(\be, w^j\al) = -\la w^j\al, \be\ra_w\varep_w(w^j\al, \be) = -\varep_w(w^j\al, \be).
\end{equation*} 
Since $\pi(\wt\al) = \pi(w^i\wt\al)$ for all $i$, we see that (\ref{eqn-2terms}) is zero as claimed.

Case 3: $(w^i\al, \be) \in \{0, \pm 1\}$ for all $i$.

Note that if $(w^i\al, \be) = 0$ for all $i$, then $\la\al, \be\ra_w = 1$, and so $\wt\al, \wt\be$ commute. Thus both sides of (\ref{eqn-rep-lift}) are zero.
Otherwise (considering the equation $(\al + w\al + w^2\al + ... + w^{d-1}\al, \be) = 0$) there exists $i$ such that $(w^i\al, \be) = -1$, and again we can assume $i = 0$. Thus the left-hand side of (\ref{eqn-rep-lift}) is 
\begin{eqnarray*}
\sum_{j \in I_{\al, \be}(-1)} \varep_w(w^j\al, \be)\rho(\pi(w^j\wt\al\wt\be))
&=& \sum_{j \in I_{\al, \be}(-1)} \varep_w(w^j\al, \be)\rho(\pi(\wt\al\wt\be)),
\end{eqnarray*}
the right-hand side is 
\begin{eqnarray*}
\rho(\pi(\wt\al)\pi(\wt\be)) - \rho(\pi(\wt\be)\pi(\wt\al)) &=& \rho(\pi(\wt\al)\pi(\wt\be)) - \la\be, \al\ra_w \rho(\pi(\wt\al)\pi(\wt\be))\\
&=& (1 - \la\be,\al\ra_w)\rho(\pi(\wt\al\wt\be)),
\end{eqnarray*}
and the theorem follows from Proposition \ref{prop-poly-relation}.
\qed

Now suppose $K = k$ and that we are in the setting of Proposition \ref{prop-basefield}, so in particular $\Gamma = \Gal(k/k_0)$ is acting on $\Heis$. If 
$V = V_0 \otimes_{k_0} k$ for some vector space $V_0$ over $k_0$ such that
$\rho: \Heis \to \GL(V)$ is $\Gamma$-equivariant for the induced action of $\Gamma$ on $\GL(V)$, then it's not hard to see that $\wt\rho$ is also $\Gamma$-equivariant and thus gives a representation of $\frakg^\Gamma$ on $V_0$.

\begin{remark}\label{rem-embedding}
Let $H$ be the split, simply connected group over $k$ with Lie algebra $\frakh$. The automorphism $\wt{w}$ of $\frakh$ induces an order-$d$ automorphism, call it $\theta$, of $H$. Let $G = H^\theta$ be the fixed-point subgroup, which is a connected reductive group with Lie algebra $\frakg$ (\cite[Theorem 8.1]{Steinberg}). Let $G^{sc}$ be the simply connected cover of $G$. Analogues of Theorem \ref{thm-rep-ext} have been used in special cases to construct an embedding $\Heis \hookrightarrow G^{sc}$ lifting the map $\Lambda_w \to \Aut(\frakg)$ described in Remark \ref{rem-innerauts} (see \cite[Section A1]{Tho16} and \cite[Theorem 3.7]{RomanoThorne2}). It is unclear whether such embeddings exist in the general case, but it seems likely that Theorem \ref{thm-rep-ext} should allow for the construction of such embeddings in many examples. For a new example of this kind, see Section \ref{section-examples}.
\end{remark}

\section{Some Lie-theoretic applications}\label{section-lie}

\subsection{On simply laced Lie algebras and lifts of automorphisms}

The following two propositions are not new results, but we point them out as general applications of the construction above.

\begin{Prop}[\cite{Lurie}]
For every simply laced root lattice $\Lambda$, Theorem 2.2 may be used to construct a simple Lie algebra of Dynkin type $\Lambda$ over $\ZZ$.
\end{Prop}

\textit{Proof}. Let $\Lambda$ be a simply laced root lattice corresponding to an irreducible root system. When $w = -1$, we have $\Lambda_w = \Lambda/2\Lambda$ and $\la \al, \be \ra_w = (-1)^{(\al, \be)}$ for all $\al, \be \in \Lambda_w$. Given a central extension $\Heis$ with commutator pairing $\la \al, \be \ra_w$, we may take $\varep$ to be trivial to obtain an input datum $(\Lambda, w, \Heis, \varep)$. With this input datum, the construction of Section \ref{section-hconstruction} reduces to the construction of \cite[Section 3.1]{Lurie}. As already noted in \cite{Lurie}, we thus obtain a canonical construction of a simple, simply laced Lie algebra over $\ZZ$ with root lattice $\Lambda$.
\qed

The next proposition was proven in more generality in \cite{AdamsHe} using very different techniques; we state it here as an easy application of the results above.

\begin{Prop}
Let $\mathfrak{s}$ be a simple, simply laced Lie algebra over an algebraically closed field $K$ of characteristic zero with Cartan subalgebra $\frakt$ and corresponding root lattice $\Lambda$. Then every elliptic automorphism $w$ of $\Lambda$ lifts to an automorphism of $\mathfrak{s}$ of the same order as $w$. 
\end{Prop}

\textit{Proof}. Let $w \in \Aut(\Lambda)$ be an elliptic automorphism of order $d$. By the considerations of Section \ref{section-pairingsandexts}, we may extend the pair $(\Lambda, w)$ to an input datum $(\Lambda, w, \Heis, \varep)$ and may use Section \ref{section-hconstruction} to produce a Lie algebra $\frakh$ with Cartan subalgebra $\frakt = \Hom(\Lambda, K)$. By Proposition \ref{prop-extendw}, we obtain an order-$d$ automorphism of $\frakh$ that restricts to $w$ on $\frakt$. We then use the fact that there exists an isomorphism $\mathfrak{s} \to \frakh$ that restricts to the identity on $\frakt$ to obtain the result. \qed

\subsection{On non-simply laced Lie algebras}\label{section-non-simply-laced}

Here we note combining Theorem \ref{thm-liealg} with Lemma \ref{lem-functor} allows us to construct every simple, non-simply laced Lie algebra over an algebraically closed field $K$ of characteristic zero. We start with the following lemma.

\begin{Lem}\label{lem-coxtrivial}
Let $\Lambda$ be a simply laced root lattice associated to an irreducible root system. Let $c$ be a Coxeter element of the Weyl group of $\Lambda$. Then $\la \cdot, \cdot \ra_c$ is trivial.
\end{Lem}

\textit{Proof}. Let $\Lambda^0 = \{\gamma \in \Lambda \otimes \mathbb{R} \mid (\gamma, \al) \in \mathbb{Z} \text{ for all } \al \in \Lambda\}$ be the weight lattice associated to $\Lambda$. By \cite[Lemma 2.3]{ReederEll}, it suffices to show that $\Lambda \subset (1 - c)\Lambda^0$. Note that $[\Lambda^0: (1-c)\Lambda^0] = \lvert\det(1 - c)\rvert = [\Lambda^0:\Lambda]$, and that $(1 - c)\Lambda^0 \subset \Lambda$ (see \cite[VI.1 Exercise 22]{Bourbaki4-6}). Thus we have $\Lambda = (1 - c)\Lambda^0$, which proves the lemma. \qed

As a side note, it's worth pointing out that even though $\la \cdot, \cdot \ra_c$ is trivial, $\varep_w$ does not seem to reduce to a simpler form in the case when $c = w$.

Suppose $\Lambda$ is type $A_{2n + 1}, D_\ell$, or $E_6$, and let $\vartheta$ be an automorphism of the Dynkin diagram of $\Lambda$. Note that there exists a Coxeter element of the Weyl group of $\Lambda$ that commutes with $\vartheta$: for an example of such an element when $\vartheta$ is nontrivial, see the table below. Let $c$ be such an element, and let $\overline\vartheta$ be the automorphism of $\Lambda_c$ induced by $\vartheta$.
Let $h$ be the Coxeter number of $\Lambda$, and let $\zeta$ be a root of unity of order $h$ in $K$. Lemma \ref{lem-coxtrivial} implies that the direct product $\Heis_c := \la \zeta \ra \times \Lambda_c$ is a central extension with commutator pairing $\la \cdot, \cdot \ra_c$.

Define $\phi \in \Aut(\Heis_c)$ by $\phi((\zeta^i, \al)) = (\zeta^i, \overline\vartheta(\al))$. By Lemma \ref{lem-epw-functor}, $(\vartheta, \phi)$ is an isomorphism of the input datum $(\Lambda, c, \Heis_c, \varep_c)$. We may use this input datum to obtain a Lie algebra $\frakh$ over $\bQ(\zeta)$.
Define $\wt\phi \in \Aut(\frakh)$ as in Lemma \ref{lem-functor}. Then $\wt\phi$ is a pinned automorphism, and $\frakh^{\wt\phi}$ is a semisimple Lie algebra with root system as listed in the table below.

In each case, the table lists the type of the root lattice $\Lambda$, a labeling of  the Dynkin diagram for $\Lambda$ with a set of simple roots $\{\al_1, \dots \al_\ell\}$, a choice of Coxeter element $c$, the order of the automorphism $\vartheta$ of the Dynkin diagram, and the Dynkin type of the resulting fixed-point Lie algebra $\frakh^{\wt\phi}$. 
We write $w_i$ for the simple reflection corresponding to $\al_i$. We note that in the case when $\lvert \vartheta \rvert = 3$, the order of $\vartheta$ doesn't pin down a unique automorphism of the Dynkin diagram of $D_4$, but any choice of order-$3$ automorphism will work in the construction.

\begin{center}
\begin{tabular}{| c | c | c | c |c |}
\hline
& & & & \\
$\Lambda$ & Dynkin diagram &  $c$ & $\lvert \vartheta \rvert$ & $\frakh^{\wt\phi}$\\
& & & & \\
\hline
& & & & \\
$A_{2n - 1}$ &  
\begin{tikzpicture}[transform shape, scale=.6]
\node[root] (b) {}; 
\node[root] (c) [right=of b] {};
\node(d) [right=of c] {$\dots$};
\node[root] (e) [right=of d] {};
\node[root] (f) [right=of e] {};
\node [above] at (b.north) {$\alpha_1$};
\node [above] at (c.north) {$\alpha_2$};
\node [above] at (e.north) {$\alpha_{2n}$};
\node [above] at (f.north) {$\alpha_{2n + 1}$};
\draw (b) -- (c);
\draw(e) -- (f);
\end{tikzpicture} & $w_1w_\ell w_2 w_{\ell - 1}\dots w_{\frac{\ell - 1}{2}}w_{\frac{\ell + 3}{2}}w_{\frac{\ell + 1}{2}}$ & 2& $C_n$\\
& & & & \\
\hline
& & & & \\
$D_{n + 1}$ 
 & 
\begin{tikzpicture}[transform shape, scale=.6]
\node[root] (b) {}; 
\node[root] (c) [right=of b] {};
\node(d) [right=of c] {$\dots$};
\node[root] (e) [right=of d] {};
\node[root] (f) [right=of e] {};
\node[root](g)[below=of e]{};
\node [above] at (b.north) {$\alpha_1$};
\node [above] at (c.north) {$\alpha_2$};
\node [above] at (e.north) {$\alpha_{\ell - 2}$};
\node [above] at (f.north) {$\alpha_{\ell - 1}$};
\node [below] at (g.south) {$\alpha_{\ell}$};
\draw (b) -- (c);
\draw(e) -- (f);
\draw(e) -- (g);
\end{tikzpicture} & $w_1w_2 \dots w_{\ell}$ & 2 & $B_n$\\ & & & & \\
\hline
& & & & \\
$D_4$ & \begin{tikzpicture}[transform shape, scale=.6]
\node[root] (c)  {};
\node[root] (d) [right=of c] {};
\node[root] (e) [right=of d] {};
\node[root] (h) [below=of d] {};
\node [above] at (c.north) {$\alpha_1$};
\node [above] at (d.north) {$\alpha_2$};
\node [above] at (e.north) {$\alpha_3$};
\node [below] at (h.south) {$\alpha_4$};
\draw (c) -- (d) -- (e);
\draw (d) -- (h);
\end{tikzpicture} 
& $w_2w_1w_3w_4$ & 3 & $G_2$\\
& & & & \\
\hline
& & & & \\
$E_6$ & 
\begin{tikzpicture}[transform shape, scale=.6]
\node[root] (b) {}; 
\node[root] (c) [right=of b] {};
\node[root] (d) [right=of c] {};
\node[root] (e) [right=of d] {};
\node[root] (f) [right=of e] {};
\node[root] (h) [below=of d] {};
\node [above] at (b.north) {$\alpha_1$};
\node [above] at (c.north) {$\alpha_3$};
\node [above] at (d.north) {$\alpha_4$};
\node [above] at (e.north) {$\alpha_5$};
\node [above] at (f.north) {$\alpha_6$};
\node [below] at (h.south) {$\alpha_2$};
\draw (b) -- (c) -- (d) -- (e) -- (f);
\draw (d) -- (h);
\end{tikzpicture} 
 & $w_2w_4w_1w_6w_3w_5$ & 2 & $F_4$\\
\hline
\end{tabular}
\end{center}

\subsubsection{Example: A $\bQ$-form of $G_2$}

In each case described in the table above, the resulting Lie algebra $\frakh^{\wt\phi}$ is defined over $\bQ(\zeta)$, but in some cases we may define an action of $\Gamma:= \Gal(\bQ(\zeta)/\bQ)$ on $\frakh^{\wt\phi}$ and use Proposition \ref{prop-basefield} to show that $(\frakh^{\wt\phi})^\Gamma$ is a Lie algebra over $\bQ$. As an example, we show that we can do this when $\frakh^{\wt\phi}$ is of type $G_2$.

Let $\Lambda$ be type $D_4$; let $\vartheta$ be a pinned automorphism of order 3, and let $c$ be the Coxeter element defined in the table for this case. Then $c$ has order $6$, so $\zeta$ is a primitive sixth root of unity. Define $\Heis_c, \phi$, and $\wt\phi$ as above, and let $\frakh$ be the Lie algebra over $\bQ(\zeta)$ resulting from the input datum $(\Lambda, c, \Heis_c, \varep_c)$. Let $\sigma$ be the nontrivial element of $\Gamma$. Then $\Gamma$ acts on $\Lambda$ via $\sigma \cdot \lambda = w_2(\lambda)$ for all $\lambda \in \Lambda$, where $w_2$ is the simple reflection corresponding to the simple root $\alpha_2$, as defined above. We note that $\sigma c \sigma^{-1} = c^{-1}$. Define an action of $\Gamma$ on $\Heis_c$ by letting $\Gamma$ act on $\la \zeta\ra$ via its natural action on $\bQ(\zeta)$ and by letting $\Gamma$ act on $\Lambda_c$ via its action on $\Lambda$. By Lemma \ref{lem-epwgalois}, the action of $\Gamma$ satisfies the conditions of Proposition \ref{prop-basefield}. Thus defining a $\Gamma$-action on $\frakh$ as in the proof of Proposition \ref{prop-basefield}, we see that $\frakh^\Gamma$ is a Lie algebra over $\bQ$. 

\begin{Lem}
The isomorphism $\wt\phi: \frakh \to \frakh$ is equivariant for the action of $\Gamma$ on $\frakh$. 
\end{Lem}

\textit{Proof}. Note that $w_2\vartheta = \vartheta w_2$. Thus for $\alpha \in \Lambda$, we have 
\begin{eqnarray*}
a_\sigma(\wt\phi(\al^\vee)) &=& a_\sigma(\vartheta(\al)^\vee)\\
&=& (w_2\vartheta(\al))^\vee\\
&=& (\vartheta w_2(\al))^\vee\\
&=& \wt\phi(a_\sigma(\al^\vee)).
\end{eqnarray*}
We have that $\wt\Lambda \simeq \la \zeta \ra \times \Lambda$. Thinking of $\wt\Lambda$ this way, an arbitrary element of $\wt\Phi$ is of the form $(\zeta^i, \al)$ for some $\al \in \Phi$. Now consider the action of $\Gamma$ on $X_{(\zeta^i, \al)}$: we have
\begin{eqnarray*}
a_\sigma(\wt\phi(X_{(\zeta^i, \al)})) &=& a_\sigma(X_{(\zeta^i, \vartheta(\al))})\\
&=& X_{(\sigma(\zeta^i), w_2\vartheta(\al))}\\
&=& X_{(\sigma(\zeta^i), \vartheta(w_2(\al))}\\
&=& \wt\phi(a_\sigma(X_{(\zeta^i, \al)})).
\end{eqnarray*}
Since $\wt\phi$ is linear and $a_\sigma$ is $\sigma$-linear, the result follows. \qed

Thus $\wt\phi$ restricts to an isomorphism $\frakh^\Gamma \to \frakh^\Gamma$, the $\Gamma$-action on $\frakh$ restricts to an action on the fixed-point subalgebra $\frakh^{\wt\phi}$, and $(\frakh^{\wt\phi})^\Gamma = (\frakh^\Gamma)^{\wt\phi}$ is a Lie algebra over $\bQ$.

\section{Examples from a simple singularity of type $E_8$}\label{section-examples}

We finish the paper by showing how input data as defined in Section \ref{section-input} naturally arise in a collection of examples coming from a simple singularity of type $E_8$. If $\Lambda$ is a root lattice of type $E_8$, then for each $d \in \{2, 3, 5\}$, there is a unique conjugacy class in $\Aut(\Lambda)$ of elliptic elements of order $d$ (see \cite[Table 1]{ReederEll}). For each of these choices of $d$, there is a naturally occurring family of curves, and in each case we obtain a central extension of the form (\ref{eqn-centralext}). We note that in each case, these families of curves have been previously studied, but we summarize some known results here to put our construction into context. For ease of exposition, we work over an algebraically closed field. 
For much of the background about elliptic surfaces related to this section see, e.g., \cite{ShiodaE8} and \cite{SchuttShioda}. 

Let $K$ be an algebraically closed field of characteristic zero, and consider the surface $\affS_0$ over $K$ defined by $X^5 + Y^3 + Z^2 = 0$, which is a simple singularity of type $E_8$. 
We may form a semiuniversal deformation $\affS \to B:= \mathbb{A}_K^8$ of this singularity with fibers of the form 
\begin{equation}\label{eqn-deformation}
X^5 + Y^3 + Z^2 + Y(c_2X^3 + c_8X^2 + c_{14}X + c_{20}) + c_{12}X^3 + c_{18}X^2 + c_{24}X + c_{30} = 0
\end{equation}
(see, e.g., \cite[Section 2.4]{Slodowy}). 
For each $\lam \in B$, we write $\affS_\lam$ for the fiber over $\lam$, and we write $B^s$ for the locus over which $\affS_\lam$ is smooth, which is the complement of a hypersurface in $B$. For each $\lam \in B^s$, we may compactify to form an elliptic surface, call it $S_\lam$, over $\mathbb{P}_K^1$, where the morphism $S_\lam \to \mathbb{P}_K^1$ is given by the $X$ coordinate. As discussed in \cite[Section 4.4]{RomanoThorne2}, $\Pic(\affS_\lam)$ is an $E_8$-root lattice.

For $d \in \{2, 3, 5\}$, let $\zeta_d \in K$ be a root of unity of order $d$.  
 For each of these $d$, there is an order-$d$ automorphism of $\affS_0$ given by 
 \begin{eqnarray*}
 (X, Y, Z) &\mapsto& (X, Y, -Z) ~~(d = 2)\\
 (X, Y, Z) &\mapsto& (X, \zeta_3Y, Z)~~(d = 3)\\
 (X, Y, Z) &\mapsto& (\zeta_5X, Y, Z)~~(d = 5). 
 \end{eqnarray*}
For each $d$, this map extends to an automorphism, call it $\theta_d$, of $\affS$ given by the same formula, inducing a $\mu_d$-action on $\affS$. If we define a $\mu_d$-action on $B$ by $\zeta_d \cdot c_i = \zeta_d^{i}c_i$, then the map $\affS \to B$ is $\mu_d$-equivariant. The fibers over the fixed points $B^{\mu_d}$ are of the form:
\begin{equation*}
X^5 + Y^3 + Z^2 + Y(c_2X^3 + c_8X^2 + c_{14}X + c_{20}) + c_{12}X^3 + c_{18}X^2 + c_{24}X + c_{30} = 0~~(d = 2)
\end{equation*}
\begin{equation*}
X^5 + Y^3 + Z^2 + c_{12}X^3 + c_{18}X^2 + c_{24}X + c_{30} = 0 ~~(d = 3)
\end{equation*}
\begin{equation*}
X^5 + Y^3 + Z^2 + c_{20}Y + c_{30} = 0 ~~(d = 5).
\end{equation*}
If $\lam \in B^{s, \mu_d}$, then $\theta_d$ restricts to an automorphism of $\affS_\lam$.

\begin{Lem}\label{lem-picelliptic}
If $\lam \in B^{s, \mu_d}$, the automorphism of $\Pic(\affS_\lam)$ induced by $\theta_d$ is elliptic. 
\end{Lem}

\textit{Proof}. 
As in \cite[Section 4.4]{RomanoThorne2}, we identify $\Pic(\affS_\lam)$ with the orthogonal complement in $\Pic(S_\lam)$ of the $\ZZ$-module generated by the zero section and the fiber above $X = \infty$. If we think of the equation (\ref{eqn-deformation}) defining $S_\lam$ as defining an elliptic curve $E_\lam$ over the field $K(X)$, then by the same logic as in \cite[Proposition 1.3]{ShiodaE8}, we have that the Mordell--Weil group $E_\lam(K(X))$ is also isomorphic to a root lattice of type $E_8$. By \cite[Lemma 4.4]{SchuttShiodaMW}, the natural surjective map $\Div(S_\lam) \to E_\lam(K(X))$ (see, e.g., \cite[Section 6.3]{SchuttShiodaMW}) factors through $\Pic(S_\lam)$, and by considering the kernel of this map \cite[(6.4)]{SchuttShiodaMW}, we see that it induces an isomorphism of lattices
\begin{equation*}
\Pic(\affS_\lam) \hookrightarrow \Pic(S_\lam) \to E(K(X)).
\end{equation*}
Thus we may think of nontrivial elements of $\Pic(\affS_\lam)$ as pairs $(Y, Z) = (a(X), b(X))$ of elements of the field $K(X)$ satisfying equation (\ref{eqn-deformation}).

Case 1: $d = 2$. 

The automorphism induced by $\theta_2$ maps $(a(X), b(X))$ to $(a(X), -b(X))$, so $(a(X), b(X))$ is fixed if and only if $b(X) = 0$. Suppose $(a(X), b(X))$ is fixed. Then $a(X)$ satisfies
\begin{equation*}
X^5 + (a(X))^3 + a(X)(c_2X^3 + c_8X^2 + c_{14}X + c_{20}) + c_{12}X^3 + c_{18}X^2 + c_{24}X + c_{30} = 0.
\end{equation*}
For this to happen, $(a(X))^3 + a(X)(c_2X^3 + c_8X^2 + c_{14}X + c_{20})$ must be a degree-5 polynomial, call it $q(x)$. Write $a(X) = \frac{a_1(X)}{a_2(X)}$ for polynomials $a_1, a_2 \in K[X]$ that have no common zeroes. We have that 
\begin{equation}\label{eqn-1}
(a_1(X))^3 + a_1(X)(a_2(X))^2(c_2X^3 + c_8X^2 + c_{14}X + c_{20}) = q(X)a_2(X)^3,
\end{equation}
so 
\begin{equation*}
(a_1(X))^3 = (a_2(X))^2(q(X)a_2(X) - a_1(X)(c_2X^3 + c_8X^2 + c_{14}X + c_{20})),
\end{equation*}
which, since $a_1$ and $a_2$ have no common zeroes, implies that $a_2 \in K$ and $a(X)$ is itself a polynomial. Now  (\ref{eqn-1}) tells us that $\deg(a_1) \leq 1$, but this contradicts the fact that $(a(X))^3 + a(X)(c_2X^3 + c_8X^2 + c_{14}X + c_{20})$ has degree 5.

Case 2: $d = 3$. 

Similarly, if $(a(X), b(X))$ is fixed, then $a(X) = 0$ and $b(X)$ must satisfy 
\begin{equation*}
X^5 + (b(X))^2 + c_{12}X^3 + c_{18}X^2 + c_{24}X + c_{30} = 0.
\end{equation*}
Using similar logic to the above, we see that $b(X)$ is a polynomial. We have a contradiction because the degree of $b(X)^2$ is even.

Case 3: $d = 5$.

The automorphism group of $\Pic(\affS_\lam)$, i.e. the Weyl group of $E_8$, contains two conjugacy classes of elements of order 5 (corresponding to diagrams $A_4$ and $A_4^2$ in \cite[Table 11]{Carter}): one consists of elliptic elements and one does not. Let $w \in \Aut(\Pic(\affS_\lam))$ be the automorphism induced by $\theta_5$.  If $w$ is not elliptic, then by finding an explicit conjugacy-class representative, we see that $w$ must fix a root in $\Pic(\affS_\lam)$. (For example, with simple roots as in the diagram below and notation for simple refections as in Section \ref{section-non-simply-laced}, we may use that $w$ is conjugate to $w_8w_7w_6w_5$.)
\begin{center}
\begin{tikzpicture}[transform shape, scale=.8] 
\node[root] (b) {}; 
\node[root] (c) [right=of b] {};
\node[root] (d) [right=of c] {};
\node[root] (e) [right=of d] {};
\node[root] (f) [right=of e] {};
\node[root] (g) [right=of f] {};
\node[root] (h) [right=of g] {};
\node[root] (i) [below=of d] {};
\node [above] at (b.north) {$\alpha_1$};
\node [above] at (c.north) {$\alpha_3$};
\node [above] at (d.north) {$\alpha_4$};
\node [above] at (e.north) {$\alpha_5$};
\node [above] at (f.north) {$\alpha_6$};
\node [above] at (g.north) {$\alpha_7$};
\node [above] at (h.north) {$\alpha_8$};
\node [below] at (i.south) {$\alpha_2$};
\draw[thick] (b) -- (c) -- (d) -- (e) -- (f)--(g)--(h);
\draw[thick] (d) -- (i);
\end{tikzpicture} 
\end{center} 

As in \cite[Proposition 1.4]{ShiodaE8} or \cite[Lemma 3.1]{ShiodaMW2}, the roots of $E_8$ are given by pairs $(a(X), b(X))$ where $a(X)$ is a polynomial of degree at most 2 and $b(X)$ is a polynomial of degree at most 3. The automorphism $w$ maps $(a(X), b(X))$ to $(a(\zeta_5X), b(\zeta_5X))$. Thus if a root $(a(X), b(X))$ is fixed, then $a(X)$ and $b(X)$ must each be constant. But $(a(X), b(X))$ must satisfy
\begin{equation*}
X^5 + a(X)^3 + b(X)^2 + c_{20}b(X) + c_{30} = 0,
\end{equation*}
a contradiction. \qed

From now on, given $d \in \{2, 3, 5\}$, we write $w_d \in \Aut(\Pic(\affS_\lam))$ for the automorphism of Lemma \ref{lem-picelliptic}. Given $d$ and given $\lam \in B^{s, \mu_d}$, let $\affC_\lam$ denote the fixed-point curve of $\theta_d$ in $\affS_\lam$. Then $\affC_\lam$ is given by the following equation, depending on $d$:
\begin{equation*}
X^5 + Y^3 + Y(c_2X^3 + c_8X^2 + c_{14}X + c_{20}) + c_{12}X^3 + c_{18}X^2 + c_{24}X + c_{30} = 0~~(d = 2)
\end{equation*}
\begin{equation*}
X^5 + Z^2 + c_{12}X^3 + c_{18}X^2 + c_{24}X + c_{30} = 0~~(d = 3)
\end{equation*}
\begin{equation*}
Y^3 + Z^2 + c_{20}Y + c_{30} = 0~~(d = 5).
\end{equation*}
In the case $d = 2$, let $C_\lam$ be the compactification of $\affC_\lam$ as in  \cite[Lemma 4.9]{Tho13}, and in the other two cases, take the projective completion of $\affC_\lam$ to yield a smooth projective curve $C_\lam$. 
In each case, there exists an isomorphism $\Pic(\affS_\lam)_{w_d} \to \Pic^0(C_\lam)[d]$ intertwining the pairing $\la \cdot, \cdot \ra_{w_d}$ with the Weil paring (see \cite[Theorem 2.10]{Zarhin} for $d = 2$ and  \cite[Proposition 4.15]{RomanoThorne2} for $d = 3$; for $d = 5$ this follows from the fact that the groups in question are isomorphic to $(\ZZ/5\ZZ)^2$ \cite[Table 1]{ReederEll}, which has a unique non-degenerate alternating pairing $(\ZZ/5\ZZ)^2 \times (\ZZ/5\ZZ)^2 \to \la \zeta_5 \ra$ up to choice of a primitive fifth root of unity). 
We now recall that in each case, we may define a central extension of $J_{C_\lam}[d]$ whose commutator pairing is given by the Weil pairing, and we may use this central extention to form an input datum.

{Case 1: $d = 2$.}

In this case, \cite[Section 1C]{Tho16} shows that given a theta characteristic $\kappa$ on $C_\lam$, we may naturally associate a central extension 
\begin{equation*}
1 \to \mu_2 \to \Heis_\kappa \to J_{C_\lam}[2] \to 1
\end{equation*}
with commutator pairing given by the Weil pairing. The special case of the construction of Section \ref{section-construction} described in \cite{Tho16} thus gives a 2-graded Lie algebra of type $E_8$. By \cite[Theorem 2.1]{Tho16}, if $C_\lam$ is defined over a subfield of $K$, then so is the corresponding Lie algebra. We note that the $2$-Selmer sets of this family of curves over $\bQ$ were studied in \cite{RomanoThorne1} by analyzing a Lie algebra of this kind, and that the construction of \cite{Tho16} was recently used in \cite{Kulkarni} in a similar setting to study Jacobians of a larger class of genus-4 curves.

{Case 2: $d = 3$.}

 In \cite[Section 4.2]{RomanoThorne2} by considering the 3-torsion subgroup of an appropriately chosen Mumford-theta group, we obtain a central extension
\begin{equation*}
1 \to \mu_3 \to \mathscr{G} \to J_{C_\lam}[3] \to 1
\end{equation*}
with commutator pairing given by the Weil pairing. Thus we may use Section \ref{section-construction} to construct a 3-graded Lie algebra of type $E_8$. 
A graded Lie algebra of this form was used in both \cite{RomanoThorne2} and \cite{RainsSam} to study the 3-Selmer groups of smooth curves of the form $C_\lam$ over $\bQ$.
If we let $\varep$ be given by $\varep(\al, \be) = (-1)^{(\al, w\be)}\la \al, \be \ra_{w_3}$ in our input datum, then we recover the construction used in \cite{RomanoThorne2}, and in this case if $C_\lam$ is defined over a subfield of $K$, then so is the corresponding Lie algebra.

{Case 3: $d = 5$.}

In this case $C_\lam$ is an elliptic curve. We may construct a Mumford-theta group as in \cite[Section 1.6]{Cremona}, obtaining a central extension
\begin{equation*}
1 \to \mathbb{G}_m \to \mathscr{G} \to {C_\lam}[5] \to 1
\end{equation*}
with commutator pairing given by the Weil pairing. The $5$-torsion subgroup $\mathscr{G}[5]$ of $\mathscr{G}$ is a central extension 
\begin{equation*}
1 \to \mu_5 \to \mathscr{G}[5] \to {C_\lam}[5] \to 1.
\end{equation*}
We have that $(\Pic(\affS_\lam), w_5, \mathscr{G}[5](K), \varep_{w_5})$ is an input datum, with $\varep_{w_5}$ as defined in Section \ref{section-epw}. Let $\frakh$ be the resulting Lie algebra and let $\wt w_5 \in \Aut(\frakh)$ be the resulting automorphism. Then $\frakh$ is a Lie algebra of type $E_8$ with 5-grading given by $\wt w_5$.
As discussed in \cite{Gross}, this 5-grading of $E_8$ gives precisely the representation studied in \cite{BS5} to determine the average size of the 5-Selmer group for elliptic curves over $\QQ$.

Note that if $C_\lam$ is an elliptic curve over $\bQ$, we may choose a finite extension $k/\bQ$ such that ${C_\lam}[5](k) = C_\lam[5](K)$ and such that $E_\lam(k(X))$ contains all 240 roots of $E_8$, which implies $E_\lam(k(X)) = E_\lam(K(X))$. We then obtain a Lie algebra $\frakh$ over $k$. The Galois group $\Gal(k/\bQ)$ naturally acts on $E_\lam(k(X))$, and it's not hard to see that the action satisfies A1 of Section \ref{section-galois}. Thinking of $\mathscr{G}$ as a subgroup of $\GL_5$ as in \cite[Section 1.6]{Cremona}, we see that the natural action of $\Gal(k/\bQ)$ on $\mathscr{G}[5](k) \subset \GL_5(k)$ restricts to the natural action of $\Gal(k/\bQ)$ on $\la \zeta_5\ra$, yielding A2 from Section \ref{section-galois}. If the isomorphism $E(k(X))_{w_5} \to C_\lam[5](k)$ intertwining $\la \cdot, \cdot \ra_{w_5}$ with the Weil pairing is equivariant for the action of $\Gal(k/\bQ)$, then by Proposition \ref{prop-basefield}, $\frakh$ descends to a Lie algebra over $\bQ$.

Our final result for the $d = 5$ case is an application of Theorem \ref{thm-rep-ext} and a new example of the kind discussed in Remark \ref{rem-embedding}. 
The analogous result for $d = 3$ is \cite[Theorem 3.7]{RomanoThorne2}, and a similar result holds for $d = 2$ by \cite[Proposition A.2]{Tho16}. 
Let $\frakh$ be the Lie algebra resulting from the input datum $(\Pic(\affS_\lam), w_5, \mathscr{G}[5](K), \varep_{w_5})$, let $\wt{w_5} \in \Aut(\frakh)$ be the lift of $w_5$ as defined in Section \ref{section-lieaut}, and let $\frakg = \frakh^{\wt{w}_5}$.
Let $H = \Aut(\frakh)$, which is a connected reductive group of type $E_8$ with maximal torus $T$ as defined in Section \ref{section-hconstruction}. The automorphism $\wt w_5 \in H$ induces an inner automorphism of $H$; denote the fixed-point subgroup by $G := H^{\wt w_5}$. Then $G$ is a semisimple connective reductive group (\cite[Theorem 8.1]{Steinberg}, \cite[Table 18]{RLYG}). Let $G^{sc}$ denote the simply connected cover of $G$. 

Let $\iota: {C_\lam}[5] \to H$ be the homomorphism of Remark \ref{rem-innerauts}, where we are identifying ${C_\lam}[5]$ with $\Pic(\affS_\lam)_{w_5}$. Using the explicit description of $\iota$ in Remark \ref{rem-innerauts}, it's not hard to check that the image of $\iota$ commutes with $\wt{w}_5$, and thus lies in $G$. Since $\la \cdot, \cdot \ra_{w_5}$ is nondegenerate \cite[Lemma 2.3]{ReederEll}, $\iota$ is injective.

\begin{Prop}\label{prop-E8(5)}
There exists a commutative diagram 
\[\begin{tikzcd}        
    1 \arrow{r} & \mu_5 \arrow{r} & G^{sc} \arrow{r} & G \arrow{r} & 1\\   
  1 \arrow{r} & \mu_5 \arrow{r}\arrow{u} & \mathscr{G}[5] \arrow{r}\arrow[hook]{u} & {C_\lam}[5] \arrow{r}\arrow[hook, "\iota"]{u} & 1
\end{tikzcd}\]
with exact rows, where $\iota$ is as defined in Remark \ref{rem-innerauts}. 
\end{Prop}

\textit{Proof}. 
The proof follows that of \cite[Theorem 3.7]{RomanoThorne2}.
By \cite[Table 18]{RLYG}, we can check that $G$ is type $A_4 \times A_4$, and by the logic in \cite[Section 3.7]{ReederTA}, we see that the center of $G$ has order 5. Thus the kernel of the natural projection $G^{sc} \to G$ is isomorphic to $\mu_5$. Let $\Heis$ denote the preimage of $\iota(\Lambda_w)$ under the projection $G^{sc} \to G$. We have a commutative diagram
\[\begin{tikzcd}        
    1 \arrow{r} & \mu_5 \arrow{r} & G^{sc} \arrow{r} & G \arrow{r} & 1\\   
  1 \arrow{r} & \mu_5 \arrow{r}\arrow{u} & \Heis \arrow{r}\arrow[hook]{u} & J_{C_\lam}[5] \arrow{r}\arrow[hook, "\iota"]{u} & 1.
\end{tikzcd}\]
Thus it suffices to show that there is an isomorphism $\mathscr{G}[5] \to \Heis$ of extensions of ${C_\lam}[5]$ by $\mu_5$. 
Note that $\mathscr{G}$ is defined as a subgroup of $\GL_5$, and $\mathscr{G}[5]$ fits into a commutative diagram 
\[\begin{tikzcd}        
    1 \arrow{r} & \mathbb{G}_m \arrow{r} & \GL_5 \arrow{r} & \PGL_5 \arrow{r} & 1\\   
  1 \arrow{r} & \mu_5 \arrow{r}\arrow[hook]{u} & \mathscr{G}[5] \arrow{r}\arrow[hook, "\rho"]{u} & {C_\lam}[5] \arrow[hook, "\overline\rho"]{u}\arrow{r} & 1.
\end{tikzcd}\]
\cite[Section 1.6]{Cremona}. 
Let $\wt\rho: \frakg \to \frakgl_5$ be the representation of $\frakg$ corresponding to $\rho$, given by Theorem \ref{thm-rep-ext}.  It's easy to see that $\wt\rho$ is irreducible with image $\fraksl_5$. 
Then $\wt\rho$ corresponds to a representation $\psi: G^{sc} \to \GL_5$ on which the center of $G^{sc}$ acts by scalars. We obtain a projective representation $\overline\psi: G \to \PGL_5$.

We claim that $\overline\rho(\lam) = \overline\psi(\iota(\lam))$ for all $\lam \in {C_\lam}[5]$. Fix $\lam \in {C_\lam}[5]$. Since the adjoint representation of $\PGL_5$ is faithful, it suffices to show that $\Ad(\overline\rho(\lam))(x) = \Ad(\overline\psi(\iota(\lam)))(x)$ for all $x \in \fraksl_5$, and since $\wt\rho$ has image $\fraksl_5$, we may assume that $x = \wt\rho(Z_{\wt\al})$ for some generator $Z_{\wt\al}$ of $\frakg$. We have
\begin{eqnarray*}
\Ad(\overline\rho(\lam))(\wt\rho(Z_{\wt\al})) &=&
\rho(\lam)\rho(\pi(\wt\al))\rho(\lam)^{-1}\\
 &=& \rho(\lam\pi(\wt\al)\lam^{-1})\\
&=& \rho(\la \lam, \al\ra_{w_5} \pi(\wt\al))\\
&=& \la \lam, \al\ra_{w_5} \wt\rho(Z_{\wt\al}),
\end{eqnarray*}
and
\begin{eqnarray*}
\Ad(\overline\psi(\iota(\lam)))(\wt\rho(Z_{\wt\al})) &=& 
\psi(\iota(\lam))\wt\rho(Z_{\wt\al})\psi(\iota(\lam))^{-1}\\
&=& \wt\rho(\iota(\lam)Z_{\wt\al}~\iota(\lam)^{-1}),
\end{eqnarray*}
which is also $\la \lam, \al\ra \wt\rho(Z_{\wt\al})$ (see Remark \ref{rem-innerauts}).

To complete the proof, it suffices to show that $\rho(\mathscr{G}[5]) = \psi(\Heis)$, since then $\psi\mid_{\Heis}$ is injective and the map $(\psi\mid_{\Heis})^{-1} \circ \rho: \mathscr{G}[5] \to \Heis$ gives the desired isomorphism of central extensions. 
Clearly $\psi(\Heis)$ is in the preimage in $\GL_5$ of $\overline\rho({C_\lam}[5])$, and $\psi(\Heis)$ is 5-torsion, so  $\psi(\Heis) \subset \rho(\mathscr{G}[5])$. The image of $\psi(\Heis)$ in $\PGL_5$ is $\overline\rho(C_\lam[5])$, so the kernel of $\psi\mid_\Heis$ is contained in $\mu_5$. So it suffices to show that $\mu_5 \subset \psi(\Heis)$. But if not then $\psi(\Heis) \simeq J_{C_\lam}[5]$ and $\psi(\Heis)$ is a maximal abelian subgroup of $\rho(\mathscr{G}[5])$, so $\mu_5 = \rho(\mu_5) \subset \psi(\Heis)$, a contradiction. 
\qed

\bibliographystyle{plain}
\bibliography{latticeexts-april2019}

\end{document}